\documentclass[amscd,amssymb,verbatim,12pt]{amsart}
\numberwithin{equation}{section}
\usepackage[normalem]{ulem}
\usepackage{epsfig,cancel}
\usepackage{multido}
\usepackage{color}
\usepackage{pstricks}
\usepackage{pst-all}
\usepackage{pst-math}
\usepackage{pst-func}
\usepackage{pst-3dplot}
\usepackage{graphicx}
\usepackage{amsmath}
\usepackage{a4,amssymb}

\newtheorem{thm}{Theorem}[section]
\newtheorem{cor}[thm]{Corollary}
\newtheorem{lem}[thm]{Lemma}

\theoremstyle{definition}

\newtheorem{rem}[thm]{Remark}

\newif\ifShowLabels
\ShowLabelstrue
\newdimen\theight
\def\TeXref#1{
     \leavevmode\vadjust{\setbox0=\hbox{{\tt
            \quad\quad  {\small  \bf #1}}}%
     \theight=\ht0
     \advance\theight  by  \dp0
     \advance\theight  by  \lineskip
     \kern -\theight \vbox  to
     \theight{\rightline{\rlap{\box0}}%
      \vss}%
      }}%

    {\begin{thm}\label{#1} \ifShowLabels \TeXref{#1} \fi}%
    {\end{thm}}

    {\begin{def}\label{#1} \ifShowLabels \TeXref{#1} \fi}%
    {\end{def}}

    {\begin{lem}\label{#1} \ifShowLabels \TeXref{#1} \fi}%
    {\end{lem}}

    {\begin{cor}\label{#1} \ifShowLabels \TeXref{#1} \fi}%
    {\end{cor}}

\newcommand{\eqRef}[1]%
     {\ifShowLabels \TeXref{#1} \fi
      \begin{equation}\label{#1} }

\ShowLabelsfalse

\newcommand{\vsp}{\vskip 1em}
\newcommand{\vspp}{\vskip 2em}
\newcommand{\NI}{\noindent}
\newcommand{\bea}{\begin{eqnarray}}
\newcommand{\eea}{\end{eqnarray}}
\newcommand{\IR}{\mathbb{R}}
\newcommand{\bas}{\begin{align*}}
\newcommand{\eas}{\end{align*}}
\newcommand{\ba}{\begin{align}}
\newcommand{\ea}{\end{align}}

\newcommand{\be}{\begin{equation}}
\newcommand{\ee}{\end{equation}}
\newcommand{\ben}{\begin{eqnarray*}}
\newcommand{\een}{\end{eqnarray*}}
\newcommand{\lam}{\lambda}
\newcommand{\Om}{\Omega}
\newcommand{\om}{\omega}
\newcommand{\tht}{\theta}

\newcommand{\al}{\alpha}
\newcommand{\bt}{\beta}
\newcommand{\bb}{\bar{\beta}}

\newcommand{\Lmn}{\Lambda_{\min}}

\newcommand{\Lmx}{\Lambda_{\max}}

\newcommand{\g}{\gamma}
\newcommand{\ve}{\varepsilon}
\newcommand{\dl}{\delta}

\newcommand{\s}{\sigma}

\newcommand{\kp}{\kappa}

\newcommand{\gs}{\gamma^*}

\newcommand{\wb}{\bar{w}}

\newcommand{\cH}{\mathcal{H}}

\title[Phragm\'en-Lindel\"of]{A Phragm\'en-Lindel\"of property of viscosity solutions to a class of nonlinear, possibly degenerate, parabolic equations}

\author[Bhattacharya and Marazzi]{Tilak Bhattacharya and Leonardo Marazzi}

\begin{document}

\maketitle

\begin{abstract} We study Phragm\'en-Lindel\"of properties of viscosity solutions to a class of doubly nonlinear 
parabolic equations in $\IR^n\times (0,T)$. We also include an application to some doubly nonlinear equations.
\end{abstract}

\section{\bf Introduction}

In this work, we discuss Phragm\'en-Lindel\"of type results for a class of nonlinear parabolic equations. This is a follow-up of the work in \cite{BM3} where we 
stated similar results for viscosity solutions of Trudinger's equation in $\IR^n\times (0,T)$, where $n\ge 2$ and $0<T<\infty$. 

We introduce notations for our discussion. Let $n\ge 2$, $g:\IR^n\rightarrow (0,\infty)$ and $h:\IR^n\rightarrow \IR$ be two continuous functions. We impose that 
\eqRef{sec1.0} 
\max\left( \sup_x |\log g(x)|,\; \sup_x |h(x)| \right)<\infty
\ee
Let $0<T<\infty$ and define $\IR^n_T=\IR^n\times (0,T)$. 

Our motivation for the work arises from the study of viscosity solutions of doubly nonlinear equations of the kind 
\eqRef{sec1.1}
H(Du, D^2u)-f(u)u_t=0, \;\mbox{in $\IR^n_T,\;\;u(x,t)>0$ and $u(x,0)=g(x)$, $\forall x$ in $\IR^n$,}
\ee
where $H$ satisfies certain homogeneity conditions and $f:\IR^+\rightarrow \IR^+$ is a non-decreasing continuous function, see Section 2 for more details. As shown in \cite{BM5}, if $f$ satisfies certain conditions then a change of variable $u=\phi(v)$ transforms (\ref{sec1.1}) to 
\eqRef{sec1.2}
H(Dv, D^2v+Z(v)Dv\otimes Dv)-v_t=0, \;\;\mbox{in $\IR^n_T$, and $v(x,0)=\phi^{-1}(g(x))$, $\forall x$ in $\IR^n$,}
\ee
where $Z:\IR\rightarrow \IR^+$ is a non-increasing function. As observed in \cite{BM2, BM5}, one can conclude a comparison principle for (\ref{sec1.2}), and hence, for (\ref{sec1.1}). 

An example of such an equation is the well-known Trudinger's equation \cite{TR}: 
\ben
\mbox{div}\left(|Du|^{p-2}Du \right)-(p-1)u^{p-2}u_t=0,\;\;\mbox{in $\IR^n_T$, and $u>0$.}
\een
The works in \cite{BM2, BM4} address the existence and uniqueness of viscosity solutions $u$, for $p\ge 2$, in cylindrical domains $\Om\times(0,T)$, where $\Om\subset \IR^n$ is a bounded domain, and \cite{BM3} includes Phragm\'en-Lindel\"of type results.

A related but some what more general equation is to consider, in $\IR^n_T$, 
\ben
&&\mbox{div}\left(|Du|^{p-2}Du \right)+\chi(t)|Du|^\s-(p-1)u^{p-2}u_t=0,\\
&& \qquad \mbox{and $u>0$, $u(x,0)=g(x)$, $\forall x$ in $\IR^n$,}
\een
where $\s\ge 0$ and $\chi(t)$ is continuous on $[0,T]$. Employing the change of variables $u=e^v$ (see \cite{BM2}), we obtain the equation 
\ben
&&\mbox{div}\left(|Dv|^{p-2}Dv \right)+(p-1)|Dv|^p+\chi(t) e^{(\s-(p-1))v} |Dv|^\s-(p-1)v_t=0,\;\;\mbox{in $\IR^n_T$,}\\
&&\qquad\mbox{and $v(x,0)=\log g(x)$, $\forall x$ in $\IR^n$.}
\een
Writing $H(Dw, D^2w)=$div$(|Dw|^{p-2}|Dw|)$, the above equation may be written as
\ben
&&H(Dv, D^2v+Dv\otimes Dv)+\chi(t) e^{(\s-(p-1))v} |Dv|^\s-(p-1)v_t=0,\;\mbox{in $\IR^n_T$,}\\
&&\qquad\mbox{$v(x,0)=\log g(x)$, $\forall x$ in $\IR^n$.}
\een
At this time, it is not clear to us as to how to address the above equation. Nonetheless, the above discussion provides motivation for addressing the following related question of studying Phragm\'en-Lindel\"of results for equations of the kind
\bea\label{sec1.3}
&&H(Dv, D^2v+Z(v)Dv\otimes Dv)+\chi(t) |Dv|^{\s}-v_t=0, \nonumber\\
&&\qquad\mbox{$v(x,0)=h(x)$, for all $x$ in $\IR^n$.}
\eea
Here $\chi,\;h$ can have any sign. 

We will show that if $v$ satisfies certain growth conditions, for large $|x|$, then $v$ satisfies a maximum principle. A similar conclusion follows for the equation in (\ref{sec1.1}). We assume $\inf_{\IR}Z(s)>0$ for the main results and this strongly influences our work. It is clear that $Z(v)Dv\otimes Dv$ and $\chi(t)|Dv|^\s$ are dueling terms and
the analysis will bear this out. Moreover, it will also show how the imposed growth rates and solutions are influenced by the power $\s$.  

We do not address existence and uniqueness issues in this work. It would be interesting to know if the growth rates stated in this work would imply such results. Omitted also from this work is the question of optimality of the growth rates.

We have divided our work as follows. In Section 2, we present some notations, assumptions and main results. In Sections 3 and 4, we present comparison principles, a change of variables result and calculations for some of the auxiliary functions we use. Sections 5 and 6 address the super-solutions and sub-solutions respectively. Finally, Section 7 presents proofs of the main results.

For additional discussion and motivation, we direct the reader to the works \cite{AJK,CIL, ED, JL}. 

\section {Notations, assumptions and main results}

We state that through out this work sub-solutions or super-solutions or solutions are understood in the viscosity sense, see \cite{BM5,CIL} for definitions. We $usc(lsc)$ for upper(lower) semicontinuous functions.

We introduce notations that will be used throughout this work. We take $n\ge 2$.
Let $0<T<\infty$ and set $\IR^n_T=\IR^n\times (0,T)=\{(x,t):\;x\in \IR^n\;\mbox{and}\;0<t<T\}$. The functions $g$ and $h$ will always satisfy (\ref{sec1.0})

By $o$, we denote the origin in $\IR^n$ and $e$ denotes a unit vector in $\IR^n$. The letters $x,\;y$ will denote points in $\IR^n$.  
Let $S^{n\times n}$ be the set of all symmetric $n\times n$ real matrices, $I$ be the $n\times n$ identity matrix and $O$ the $n\times n$ zero matrix. 

We now describe the conditions placed on $H$. 
\vsp
{\bf Condition A (Monotonicity):} The operator $H:\IR^n\times S^n\rightarrow \IR$ is continuous for any $(q,X)\in \IR^n\times S^{n\times n}$. We assume that
\bea\label{sec2.1}
&&\mbox{(i)}\;H(q,X)\le H(q,Y),\;\mbox{for any $q\in \IR^n$ and for any $X,\;Y$ in $S^{n\times n}$ with $X\le Y$}, \nonumber\\
&&\mbox{(ii)}\;H(q,O)=0,\;\mbox{for any $q\in \IR^n$.}
\eea
Clearly, for any $q\in \IR^n$ and $X\in S^{n\times n}$, $H(q,X)\ge 0$ if $X\ge O$.
\vsp
{\bf Condition B (Homogeneity):} There is a constant $k_1\ge 0$ such that for any $(q,X)\in \IR^n\times S^{n\times n}$,
\bea\label{sec2.2}
&&\mbox{(i)}\;H(\tht q, X)=|\tht|^{k_1}H(q, X),\;\mbox{for any $\tht\in \IR$, and}\nonumber\\
&&\mbox{(ii)}\;H(q, \tht X)=\tht H(q, X),\;\mbox{for any $\tht>0.$}
\eea
Our results in this work can be adapted to include the case $H(q, \tht X)=\tht^{k_2}H(q, X)$ where $k_2$ is an odd natural number. However, in this work, $k_2=1$. We note that if $k_1=0$ then
$H(q, X)=H(q/\tht, X),\forall \tht>0$. Hence, $H(q, X)=H(X).$

\vsp
Before stating the next condition, we introduce additional notation. Let $\rho\in \IR^n$ be a vector and we write its component form as $(\rho_1,\rho_2,\ldots, \rho_n)$. Recall that $(\rho\otimes \rho)_{ij}=\rho_i\rho_j,\;i,j=1,\ldots n$. Clearly, $\rho\otimes \rho\in S^{n\times n}$ and $\rho\otimes \rho\ge O$.

Recalling that $e\in \IR^n$ is a unit vector, define, for every $\lam\in \IR,$
\bea\label{sec2.4}
\Lambda_{\min}(\lam)=\min_{e}H(e, \lam e\otimes e-I)\;\;\mbox{and}\;\;\Lambda_{\max}(\lam)=\max_{e}H(e, \lam e\otimes e+ I). 
\eea
By Condition A, $\Lmn(\lam)$ and $\Lmx(\lam)$ are both non decreasing functions of $\lam$. 
\vsp
{\bf Condition C(Growth at Infinity):} Firstly, we require that
$$\max_{e}H(e, -I)<0<\min_{e}H(e, I).$$
Next, we assume that $H$ satisfies 
\bea\label{sec2.5}
\Lambda_{\min}(\lam_0)=\min_{e}H(e, \lam_0 e\otimes e-I)>0,\;\mbox{for some $\lam_0>1$}. 
\eea
We require $\lam_0>1$ since, by Condition A, $e\otimes e-I \le O$.
\vsp
We state some simple implications of Condition C. 
By Condition A, $\Lambda_{\min}(\lam)\ge\Lambda_{\min}(\lam_0)>0,$ for any $\lam\ge \lam_0$.
By Condition B, for $\lam\ge \lam_0$,
$$\Lambda_{\min}(\lam)= \left( \frac{\lam}{\lam_0} \right) \min_{e}H \left(e, \lam_0 e\otimes e-\frac{\lam_0}{\lam} I \right)\ge \frac{\lam \Lambda_{\min}(\lam_0)}{\lam_0},$$
since $\lam_0e\otimes e-I\le \lam_0e\otimes e-(\lam_0/\lam) I$.
Clearly, 
\eqRef{sec2.50}
\frac{\Lambda_{\min}(\lam) }{\lam}\ge \frac{\Lambda_{\min}(\lam_0)}{\lam_0} \quad\mbox{and}\quad \sup_{\lam>0}\Lambda_{\min}(\lam)=\infty.
 \ee
Thus, under Conditions A, B and C, (\ref{sec2.5}) implies (\ref{sec2.50}). Clearly, (\ref{sec2.50}) implies (\ref{sec2.5}). 

Next, by Conditions A, B and (\ref{sec2.50}), for $\lam\ge \lam_0$,
\eqRef{sec2.51}
\min_eH(e, e\otimes e)\ge \frac{\Lambda_{\min}(\lam)}{\lam}=\min_{e}H \left(e, e\otimes e - \frac{I}{\lam}  \right)\ge \frac{\Lambda_{\min}(\lam_0)}{\lam_0}>0.
\ee
If $\min_e H(e, e\otimes e)>0$ then by the continuity of $H$, Conditions A and B, $\min_eH(e, \lam_0 e\otimes e-I)>0$ for some $\lam_0>1$. See Section 3 for further discussion. 
\vsp
Examples of operators that satisfy Conditions A, B and C are the $p$-Laplacian, pseudo $p$-Laplacian, for $p\ge 2$, infinity-Laplacian and the Pucci operators, see \cite{BM5} for a more detailed discussion. It is easily seen that they all satisfy (\ref{sec2.51}). We remark that some of the conditions here differ from those in \cite{BM5}.  
\vsp
For the rest of this work, we set
\eqRef{sec2.6}
k=k_1+1\quad\mbox{and}\quad \g=k_1+2=k+1.
\ee
Also, $\chi:[0,T]\rightarrow \IR$ is a continuous function and $Z:\IR\rightarrow \IR^+$ is a non-increasing continuous function with $0<\inf Z\le \sup Z<\infty$. Let $h:\IR^n\rightarrow \IR$, continuous, satisfy (\ref{sec1.0}).
\vsp
We now state the main results of this work. 
For Theorems \ref{sec2.10} and \ref{sec2.11}, we assume that Conditions A, B and C hold. 
We set
\eqRef{sec2.61}
\mu=\inf_{x\in \IR^n} h(x),\;\; \nu=\sup_{x\in \IR^n} h(x),\;\;\ell=\inf_sZ(s),\;\; \cH=\min_{|e|=1}H(e, e\otimes e)\;\;\mbox{and}\;\;\al=\sup_{0\le t\le T} |\chi(t)|. \ee
\vsp
\begin{thm}\label{sec2.10}{(Maximum Principle)} Let  $0<T<\infty$, and $\nu$ and $\al$ be as in (\ref{sec2.61}).
Assume that $\sup_{\IR^n} \nu<\infty$.
Let $u\in usc(\IR^n_T)$ solve
\ben
&&H(Du, D^2u+Z(u)Du\otimes Du)+\chi(t)|Du|^\s-u_t\ge 0,\quad\mbox{in $\IR^n_T$}\\
&&\mbox{and $u(x)\le h(x),\;\forall x\in \IR^n$}.
\een 
Suppose that there is $\dl>0$ such that   
\ben
\sup_{0\le |x|\le R,\;0\le t\le T}u(x,t)\le   o( R^\dl),\;\;\mbox{as $R\rightarrow \infty$.}
\een
The following hold.

(a) Let $\s=0$. Either (i) $k=1$ i.e., $\g=2$ and $\dl=2-\ve$, for any fixed and small $\ve>0$, or (ii) $k>1$ and $\dl=\g/k$. In both cases,
$$\sup_{\IR^n_T} u(x,t)\le \nu+\al t.$$

(b) Let $0<\s\le \g$.  Either (i) $k=1$ i.e., $\g=2$ and $\bt=2-\ve$, for any fixed and small $\ve>0$, or (ii) $k>1$ and $\bt=\g/k$.  In both cases,
$$\sup_{\IR^n_T} u(x,t)\le \nu.$$

(c) Let $\s>\g$ and $\dl=\s/(\s-1)$. Then 
$$\sup_{\IR^n_T} u(x,t)\le \nu.\;\;\;\;\;\;\Box$$
\end{thm}
\vsp
\begin{thm}\label{sec2.11}{(Minimum Principle)} Let  $0<T<\infty$ and $\mu, \;\al,\;\ell$ and $\cH$ be as in (\ref{sec2.61}).
Assume that $\mu>-\infty$.
Let $u\in lsc(\IR^n_T)$ solve
\ben
&&H(Du, D^2u+Z(u)Du\otimes Du)+\chi(t)|Du|^\s-u_t\le 0,\quad\mbox{in $\IR^n_T$,}\\
&&\mbox{and $u(x)\ge h(x),\;\forall x\in \IR^n$}.
\een 

(a) If $\s=0$ then $\inf_{\IR^n_T} u(x,t)\ge \mu-\al t.$
\vsp
(b) If $0<\s<\g$ then $\displaystyle{\inf_{\IR^n_T} u(x,t)\ge \mu- \left( \al^\g/(\ell \cH)^\s \right)^{1/(\g-\s)} t.}$ 
\vsp
If $\chi(t)\ge 0$ then, for any $0\le \s<\infty$, $\inf_{\IR^n_T}u(x,t)\ge \mu.$ 
\vsp
(c) If $\s=\g$ and $\al< \ell \cH$ then $\inf_{\IR^n_T} u(x,t)\ge \mu$. 
\vsp
(d) Let $\s=\g$ and $\al\ge \ell \cH$. Assume that either (i) 
$k=1$ ($\g=2$) and, for any fixed small $\ve>0$, we have
\ben
\sup_{0\le |x|\le R,\;0\le t\le T}(-u(x,t)) \le   o( R^{2-\ve}),\;\;\mbox{as $R\rightarrow \infty$,},
\een
or (ii) $k>1$ ($\g>2$) and 
\ben
\sup_{0\le |x|\le R,\;0\le t\le T}(-u(x,t)) \le   o( R^{\g/k}),\;\;\mbox{as $R\rightarrow \infty.$}
\een
Then $\inf_{\IR^n_T} u(x,t)\ge \mu.$
\vsp
(e) If $\s>\g$ and  
\ben
\sup_{0\le |x|\le R,\;0\le t\le T}(-u(x,t)) \le   o( R^{\s/(\s-1)}),\;\;\mbox{as $R\rightarrow \infty$},
\een
then $\inf_{\IR^n_T} u(x,t)\ge \mu.$
\end{thm}

As an observation, if $H$ is quasilinear (the $p$-Laplacian, for instance) then
$$H(Dw, D^2w+Z(w)Dw\otimes Dw)=H(Dw, D^2w)+Z(w)|Dw|^\g H(e, e\otimes e).$$
The above holds for both $Dw\ne 0$ and $Dw=0$ since $H(q, O)=0,$ for any $q\in \IR^n$. 
If $H$ is the $p$-Laplacian, $p\ge 2$, $\chi(t)=1-p$ and $Z(w)=1$ then $k_1=p-2$,
$k=p-1$, $\g=p$ and $H(e, e\otimes e)=p-1.$ Clearly,
\ben
&&H(Dw, D^2w+Dw\otimes Dw)-(p-1)|Dw|^p\\
&&\quad\qquad\qquad\qquad=H(Dw, D^2w)+(p-1)|Dw|^p-(p-1)|Dw|^p=H(Dw, D^2w).
\een
Thus, the above results also apply to equations of the kind $\Delta_pu-u_t=0$ in $\IR^n_T$.
\vsp
Finally, we obtain the following theorem for a class of doubly nonlinear equations. We apply parts (a) of Theorems \ref{sec2.10} and \ref{sec2.11} with $\al=\s=0$.

\begin{thm}\label{sec2.12} Let $k\ge 1$, $f: [0,\infty)\rightarrow [0,\infty)$ be a $C^1$ non-decreasing function, and $g:\IR^n\rightarrow (0,\infty)$ be such that $0<\inf_x g(x)\le \sup_x g(x)<\infty.$

{ Let $k>1$. We assume that $f^{1/(k-1)}$ is concave and 
$$0<\inf_{0\le s<\infty} \frac{d}{ds} f^{1/(k-1)}(s)\le \sup_{0\le s<\infty} \frac{d}{ds} f^{1/(k-1)}(s)<\infty.$$}
Select
$\phi:\IR\rightarrow [0,\infty)$, a $C^2$ increasing function, such that 
$$\phi^{\prime}(\tau)= f(\phi(\tau))^{1/(k-1)}. $$

(a) Suppose that $u\in usc( \overline{\IR^n_T}),\;u>0,$ solves
$$H(Du, D^2u)-f(u) u_t\ge 0,\;\mbox{in $\IR^n_T$ and $u(x,0)\le g(x),\;\forall x\in \IR^n$.}$$
Suppose that $\sup_{|x|\le R, \;0\le t\le T} u(x,t)\le \phi( o(R^{\g/k}))$, as $R\rightarrow \infty$. Then
$\sup_{ \IR^n_T}u(x,t)\le \sup_x g(x).$

(b) Suppose that $u\in lsc( \overline{\IR^n_T}),\;u>0,$ solves
$$H(Du, D^2u)-f(u) u_t\le 0,\;\mbox{in $\IR^n_T$ and $u(x,0)\ge g(x),\;\forall x\in \IR^n$.}$$
Then
$\inf_{ \IR^n_T}u(x,t)\ge \inf_x g(x).$

If $k=1$, we take $f\equiv 1$ and $\phi(\tau)=e^\tau$. The conclusion in part (a) holds provided that we assume that, for any $\ve>0$, $\sup_{|x|\le R, \;0\le t\le T} u(x,t)\le \exp( o(R^{2-\ve}) )$, as $R\rightarrow \infty$. The conclusion in part (b) holds without any modifications.

\end{thm}
{ The condition placed on $f^{1/(k-1)}$ implies that $\phi^{\prime\prime}(\tau)/\phi^{\prime}(\tau)$ is positive and non-increasing in $\tau\in (-\infty, \infty)$. Moreover, this quotient is bounded from above and its lower bound is positive. See Section 3.}

\section{Preliminaries}

In this section, we present some calculations important for our work, a comparison principle and a change of variable result useful for our work. We also present additional discussion about the condition in (\ref{sec2.5}).

For definitions and a discussion of viscosity solutions, we direct the reader to \cite{CIL} and Section 2 in \cite{BM3}. For additional discussion and motivation, see \cite{AJK, BM2, ED, JL}.

Recall that  $Z:\IR\rightarrow \IR^+$ is continuous and non-increasing. We assume  that
\eqRef{sec3.0}
0<\inf_{\IR} Z\le \sup_{\IR} Z<\infty.
\ee
\vsp
We present now some elementary but important calculations.
Let $z\in \IR^n$ and $r=|x-z|$. Suppose that $v(x)=v(r)$ is a $C^2$ function. Set $e=(e_1,e_2,\cdots,e_n)$ where $e_i=(x-z)_i/r,\;\forall i=1,2,\cdots, n$. For $x\ne z$, 
\eqRef{sec3.1}
\left\{\begin{array} {lcr} Dv=v^{\prime}(r)e,  \;\;\;\; Dv\otimes Dv=(v^{\prime}(r))^2e\otimes e,\quad\mbox{and}\\
D^2v=(v^\prime(r)/r)\left(I-e\otimes e\right)+v^{\prime\prime}(r) e\otimes e.\end{array}\right.
\ee
\vsp
\begin{rem}\label{sec3.2} Let $\kappa:[0,T]\rightarrow (0,\infty)$ be a $C^1$ function and $Z$ be as in (\ref{sec3.0}). Fix $z\in \IR^n$ and set $r=|x-z|$. Suppose that $0<\mathcal{R}\le \infty$ and 
$w:B_{\mathcal{R}}(z)\times [0,T]\rightarrow \IR$ is a $C^1$ function with $w(x,t)=w(r,t).$ Assume also that $w$ is $C^2$ in $x$ in
$B_{\mathcal{R}}(z)\setminus \{z\}$.

Using (\ref{sec3.1}) in $x\ne z$, we get that
\bea\label{sec3.3}
&&H(Dw, D^2w+Z(w) Dw\otimes Dw)\nonumber\\
&&\quad\qquad\qquad\qquad\qquad\qquad=H\left( w_r e, \;\frac{w_r}{r}\left( I-e\otimes e\right) +\left(w_{rr} +Z(w)(w_r)^2 \right)e\otimes e \; \right).
\eea
\vsp
Recall Condition B in (\ref{sec2.2}) and (\ref{sec2.6}) i.e, $k=k_1+1$ and $\g=k_1+2$. 

{\bf Case (a) $w_r> 0$:}  Let $a$ be any scalar, $b\ge0$ and  $\mathcal{R}=\infty$. Suppose that $w(x,t)=(a+bv(r))\kappa(t),$ where $v^{\prime}(r)>0$ and $\kappa\ge 0$. 

In (\ref{sec3.3}), factor $w_r$ from the first entry, $w_r/r$ from the second, use (\ref{sec2.2}) and $k=k_1+1$ to get
\bea\label{sec3.4}
&&H(Dw, D^2w+Z(w) Dw\otimes Dw)=\frac{ w_r^k}{ r}  H\left( e, I +\left( \frac{r w_{rr}(r) }{w_r} + r w_r Z(w)-1\right)e\otimes e\right) \nonumber\\
&&\quad\qquad\qquad=\frac{(b v^{\prime}(r)\kappa(t))^k} {r}   H\left( e, \;I +\left( \frac{rv^{\prime\prime}(r)}{v^{\prime}(r)}-1 \;+\;b\kappa(t)(rv^{\prime}(r)) Z(w) \; \right)e\otimes e\;\right).
\eea
This version will be used for small $r$.
\vsp
For the second version, in (\ref{sec3.3}) we factor $w_r$ from the first entry, $w_r^2$ from the second entry of $H$, use (\ref{sec2.2}) and $\g=k_1+2$ to get, in $r>0$,
\bea\label{sec3.5}
&&H(Dw, D^2w+Z(w) Dw\otimes Dw) \nonumber\\
&&\qquad\qquad \qquad=w_r^\g H\left( e, \;\frac{I-e\otimes e}{r w_r}+\left( \frac{w_{rr} }{w_r^2 } + Z(w) \right)e\otimes e\;\right) \nonumber\\
&&\qquad\qquad\qquad=(bv^{\prime}(r)\kappa(t) )^\g H\left( e, \; \frac{ I -e\otimes e }{ b\kappa(t)(r v^{\prime}(r)) } \;+   \left( \frac{v^{\prime\prime}(r)}{b \kappa(t)(v^\prime(r))^2 } + Z(w) \right)e\otimes e   \right).
\eea
This version will be used for large $r$.

In this work, we take $0<b<1$. By factoring $1/b$ from the second entry in $H$, using Condition B and $\g=k+1$, the above may be rewritten as
\bea\label{sec3.7}
&&H(Dw, D^2w+Z(w) Dw\otimes Dw) \nonumber\\
&&\qquad\qquad\qquad=b^k(\kappa(t)v^{\prime}(r) )^\g H\left( e, \; \frac{ I -e\otimes e }{ \kappa(t) rv^{\prime}(r) } \;+
 \left( \frac{v^{\prime\prime}(r)}{ \kappa(t)(v^\prime(r))^2 } +b Z(w) \right)e\otimes e   \right).
\eea
\vsp
{\bf Case (b) $w_r< 0$:} Using (\ref{sec2.2}), (\ref{sec3.3}) and arguing as in part (a), we get
\bea\label{sec3.8}
&&H(Dw, D^2w+Z(w) Dw\otimes Dw) \nonumber\\
&&\quad\qquad\qquad=\frac{ |w_r|^k}{ r } H\left( e, \;\left(\; r|w_r|Z(w)\;+1   +\frac{rw_{rr}}{|w_r|}\right)e\otimes e -I \right)\nonumber\\
&& \quad\qquad\qquad=|w_r|^\g H\left( e, \;\frac{I-e\otimes e}{r w_r }+ \left( \frac{ w_{rr}  }{ w_r^2 } + Z(w) \right)e\otimes e\;\right).
\eea
\vsp

Set $w(x,t)=v(r)-\kappa(t)$, where $v^{\prime}(r)< 0$. Using (\ref{sec3.8}) we get
\bea\label{sec3.80}
&&H(Dw, D^2w+Z(w) Dw\otimes Dw) \nonumber\\
&&\qquad\qquad=\frac{ |v^{\prime}(r)|^k}{ r } H\left( e, \;\left(\; r| v^{\prime}(r)|Z(w)\;+1   - \frac{r v^{\prime\prime}(r)}{v^\prime(r)}\right)e\otimes e -I \right)\nonumber\\
&& \qquad\qquad=|v^{\prime}(r)|^\g H\left( e, \;\frac{I-e\otimes e}{r v^{\prime}(r) }+ \left( \frac{v^{\prime\prime} (r)}{ (v^{\prime}(r) )^2 } + Z(w) \right)e\otimes e\;\right).\;\;\;\Box
\eea
\end{rem}
\vsp

We now state a comparison principle that will be used in this work. See \cite{CIL} and Section 4 in \cite{BM5}.
\vsp
Let $F:\IR^+\times \IR\times \IR^n\times S^{n\times n}\rightarrow\IR$ be continuous. Suppose that $F$ satisfies $\forall X,\;Y\in S^{n\times n}$, with $X\le Y$,
\bea\label{sec3.17}
\mbox{$F(t, r_1, p, X)\le F(t, r_2, p, Y),$ $\forall (t,p)\in \IR^+\times \IR^n$ and $r_1\ge r_2$.}
\eea

\begin{lem}\label{sec3.18}{(Comparison principle)} Let $F$ satisfy (\ref{sec3.17}), $g:\IR\rightarrow \IR$ be a bounded non-increasing continuous function and $\hat{f}:\IR^+\rightarrow \IR^+$ be continuous.
Suppose that $\Om\subset \IR^n$ is a bounded domain and $T>0$. Let $u\in usc(\Om_T\cup P_T)$ and $v\in lsc(\Om_T\cup P_T)$ satisfy in $\Om_T$,
\ben
&&F(t,u, Du, D^2u+g(u)Du\otimes Du)- \hat{f}(t)u_t\ge 0\;\;\mbox{and}\\
&&F(t,v, Dv, D^2v+g(v)Dv\otimes Dv)-  \hat{f}(t)v_t\le 0.
\een
If $\sup_{P_T}v<\infty$ and $u\le v$ on $P_T$ then $u\le v$ in $\Om_T$. \quad $\Box$
\end{lem}
\vsp
We now discuss a change of variables result in the context of doubly nonlinear equations of the kind: 
$$H(Du, D^2u)-f(u)u_t=0,\;\;\;\mbox{where $u>0$.}$$
This is shown in Lemma 2.3 in \cite{BM5}. An earlier version appears in \cite{BM3}. 

Recall that $k=k_1+1$ and $\g=k_1+2$.
Let $f:[0,\infty) \rightarrow [0,\infty)$ be a $C^1$ increasing function.
For $k>1$, the function $\phi:\IR\rightarrow [0,\infty)$ be a $C^2$ solution of
\bea\label{sec3.22}
\frac{d\phi}{d\tau}=\left\{(f\circ \phi)(\tau)\right\}^{1/(k-1)},\;\;\phi\ge 0.
\eea
Thus, $\phi$ is increasing. We will assume further that
\eqRef{sec3.23}
\mbox{$f^{1/(k-1)}$ is concave, i.e,}\;\left\{f^{1/(k-1)}\right\}^{\prime}{ (s)}\;\;\mbox{is non-increasing in $s$. }
\ee
For example, if $f(s)=s^\al,\;\al\ge 0,$ then $f^{1/(k-1)}$ is concave if $\al\le k-1$. This condition ensures that the comparison principle in Lemma \ref{sec3.18} holds.

Using (\ref{sec3.22}) and (\ref{sec3.23}) we get that
$$\frac{\phi^{\prime\prime}(\tau)}{\phi^{\prime}(\tau)}=\left( \frac{d}{d\phi} f(\phi)^{1/(k-1)}  \right)(\phi(\tau))\;\;\mbox{is non-increasing in $\tau$.}$$

{Note that if $f(s)=s^\al,0<\al<k-1$, then $(\phi^{\prime\prime}/\phi^\prime)(s)=C/s$, for an appropriate $C=C(\al,k)$. Our work, however, excludes such cases as the quotient becomes small for large $s$. If $\al=k-1$ then
$\phi(\tau)=Ae^\tau$ and $(\phi^{\prime\prime}/\phi^\prime)(\tau)=1$. The latter is included in our work and is addressed in Theorem \ref{sec2.12}. } 

We now state the following change of variables lemma which is a simplified version of Lemma 2.3 in \cite{BM5}.

\begin{lem}\label{sec3.24} Let $H$ satisfy Conditions A and B, see (\ref{sec2.1}) and (\ref{sec2.2}) and $f:[0,\infty)\rightarrow [0,\infty)$ be a $C^1$ increasing function that satisfies
(\ref{sec3.23}). Suppose that $\phi:\IR\rightarrow [0,\infty)$ is a positive $C^2$ increasing function. 

Case (i): Suppose that $k>1$ and $\phi$ is as in (\ref{sec3.22}). We assume that $f$ is non-constant, $u>0$ and $v=\phi^{-1}(u)$.

Then $u\in usc(lsc)(\Om_T)$ solves $H(Du, D^2u)\ge(\le) f(u) u_t$, in $\Om_T$, if and only if $v \in usc(lsc)(\Om_T)$ and
$$H\left(Dv, D^2v+\frac{\phi^{\prime\prime}(v)}{\phi^{\prime}(v)}Dv\otimes Dv \right)\ge(\le) v_t,\;\;\mbox{in $\Om_T$}.$$

Case (ii): Let $k=1$, i.e., $k_1=0$. If $f\equiv 1$ then the claims in (a) and (b) hold if $\phi(\tau)$ is any increasing positive $C^2$ function. In particular, if $\phi(\tau)=e^\tau$ and 
$u\in usc(lsc)(\Om_T)$ then $H(D^2u)\ge(\le)u_t$ if and only if $v\in usc(lsc)(\Om_T)$ and
$H(D^2v+Dv\otimes Dv)\ge(\le)v_t.$
\end{lem}
\vsp
Finally, we make further comments on the dependence of $H(e, \lam e\otimes e\pm I)$ on $\lam$. 

\begin{rem}\label{sec3.20} Recall (\ref{sec2.1}), (\ref{sec2.2}), (\ref{sec2.4}) and (\ref{sec2.6}). As observed in Section 2 (see discussion immediately following
Condition C), 
$H(e, \lam e\otimes e \pm I)$ is non-decreasing in $\lam$, and, for $\lam\ge 0$ and any $e$,
$$\lim_{\lam\rightarrow \infty} \frac{H(e, \lam e\otimes e\pm  I)}{\lam}=\lim_{\lam\rightarrow \infty}H\left(e, e\otimes e \pm \frac{I}{\lam}\right)=H(e, e\otimes e)\ge 0.$$
This follows since $e\otimes e\ge O$ and $H(e, e\otimes e)\ge 0$.

The $p$-Laplacian, the pseudo $p$-Laplacian, the infinity-Laplacian and the Pucci type operators all satisfy $\sup_\lam\Lmn(\lam)=\infty$ and $H(e, e\otimes e)>0$ (note that eigenvalues of $e\otimes e$ are $1$ and $0$ ($0$ has multiplicity $n-1$)). Our current work applies to these operators. See Section 3 in \cite{BM5}.

Note that the condition 
$\min_eH(e, e\otimes e)>0$ implies that $\sup_{\lam} \Lmn(\lam)=\infty$.
Clearly, if $\sup_{\lam}\Lmx(\lam)<\infty$ then $H(e, e\otimes e)=0.$
Moreover, if $H$ is quasilinear then 
$$H(e, \lam e\otimes e\pm I)=H(e, \lam e\otimes e)+H(e, \pm I)=H(e, \pm I),\;\forall \lam\ge 0.$$

An example of such an operator is 
$$H(p, X)=|p|^{k_1}\left( |p|^2\mbox{Trace}(X)-\sum_{i,j=1}^n p_ip_jX_{ij} \right),\;\;\;\forall (p, X)\in \IR^n\times S^{n\times n},\;\;\forall k_1\ge 0,$$
i.e.,
$H(Du, D^2u)=|Du|^{k_1+2}\Delta u-|Du|^{k_1}\Delta_\infty u.$ Clearly, $H$ is elliptic and $\forall e$, 
$$H(e, e\otimes e)=0,\;\;\mbox{and}\;\;\;H(e, \lam e\otimes e\pm I)=\pm(n-1),\;\;\forall \lam.$$
Our current work omits such operators.

Note that 
the condition $\max_eH(e, e\otimes e)=0$ does not imply the boundedness of $H(e, \lam e\otimes e\pm I)$. An example is
$H(e,X)=det(X)$. The eigenvalues of $I+\lam e\otimes e$ are $1+\lam$ and $1$, the latter has multiplicity $n-1$. Thus,
$$H(e, \lam e\otimes e+I)=1+\lam=o(\lam^n),\;\mbox{as $\lam\rightarrow \infty$.} \;\;\;\Box$$
\end{rem}

\section{\bf Auxiliary Functions}

In this section, we record observations about auxiliary functions that are used in the proofs of the theorems in this work. We recall that
$k=k_1+1$ and $\g=k_1+2=k+1$. 

\begin{lem}\label{sec3.10} Let $1<\bb<\bt$ and $y\in \IR^n$. Set $r=|x-y|$ for all $x\in \IR^n$. 
For $r\ge 0$, define
$$v(r)=\int_0^{ r^\bt } (1+\tau^p)^{-1} \;d\tau,\;\;\;\;\mbox{where}\;\;p=\frac{\bt-\bb}{\bt}. $$

Then (i) $0<p<1$, \;\;(ii) $(1-p)\bt=\bb,$ and
$$\mbox{(iii)}\;\;r^{\bt} (1+R^{\bt p} )^{-1}\le v(r)\le r^{ \bt },\;\;\mbox{$\forall 0\le r\le R$, and, for any $R>0$.} $$

\vsp
Set $c_p=[2(1-p)]^{-1}$. If $R>1$ then
\ben
(iv)\; c_p\left( r^{\bb} -R^{\bb} \right)\le v(r)-v(R)\le (2c_p)\left( r^{\bb} -R^{\bb} \right), \;\;\mbox{for} \;r\ge R.
\een
\vsp
Moreover, in $r>0$, we have
\ben
&& (v)\;\;v^{\prime}(r)=\frac{\bt r^{\bt-1}}{1+ r^{p\bt}}\le \bt\min\left( r^{ \bb-1}, r^{\bt-1}\right),\quad (vi)\;\; rv^{\prime}(r)= \frac{\bt r^{\bt} }{ 1+r^{p\bt} }\le \bt \min\left(r^{ \bb },\;r^{\bt}\right),\\
&&\mbox{and (vii)}
\;\;\;v^{\prime\prime}(r) =\bt r^{\bt-2} 
\left[ \frac{ \bt-1+(\bb-1) r^{p\bt} }{  ( 1+r^{p\bt})^2 } \right].
\een

Next, we have 
\ben
&&(viii)\;\frac{ (v^{\prime}(r))^k }{ r}= \left( \frac{ \bt }{ 1+ r^{ p\bt} }\right)^k   r^{k\bt-\g}\le \bt^k\min\left(r^{k\bt-\g},\;\;r^{k\bb -\g} \right),\;\forall r>0,\\
&&\mbox{and (ix)}\;\;\;\bb-1\le \frac{rv^{\prime\prime}(r)}{v^{\prime}(r)}= \frac{\bt-1+(\bb-1)r^{p\bt} } { 1+ r^{p\bt} }\le \bt-1.
\een
Finally,
\ben
(x)\;\frac{v^{\prime\prime}(r)}{(v^\prime(r))^2}= \frac{\bt-1 }{\bt r^\bt } + \frac{ \bb-1 }{ \bt r^{ \bb } },\;\mbox{$\forall r>0$, and (xi)}\;
\;\frac{\bb-1}{\bt r^{\bb}}\le \frac{v^{\prime\prime}(r)}{(v^\prime(r))^2}\le \frac{2(\bt-1)}{\bt r^{\bb} },\;\forall r\ge 1.
\een
\end{lem}
\vsp
\begin{proof}
Parts (i) and (ii) follow easily. Part (iii) is a consequence of the bound $1+\tau^p\le 1+R^{p\bt},\;\forall \tau\le R^\bt.$ 
Part (iv) follows by noting that $\tau^p\le 1+\tau^p\le 2\tau^p,\;\tau\ge 1$, (ii) and writing
$$v(r)=v(R)+\int_{R^\bt}^{r^\bt} (1+\tau^p)^{-1} d\tau.$$

Parts (v), (vi) and (viii) are easily obtained by the estimate $1+r^{p\bt}\ge \max(1,\;r^{p\bt})$ and noting that $\g=k+1$ and $\bt-\bb=p\bt.$ 

To see (vii), we differentiate (v) and use (ii) to find
\ben
&&v^{\prime\prime}(r)=\bt\left[ \frac{(\bt-1) r^{\bt-2} }{ 1+r^{p\bt} }-\frac{ p\bt r^{p\bt+\bt-2} }{ (1+r^{p\bt} )^2} \right]=\bt r^{\bt-2} \left[ \frac{(\bt-1) (1+r^{p\bt})-p\bt r^{p\bt } }{ (1+r^{p\bt} )^2} \right] \\
&&\qquad\;\;=\bt r^{\bt-2} \left[ \frac{\bt-1+(\bb-1) r^{p\bt } }{ (1+r^{p\bt} )^2} \right].
\een

Applying (v), (vii) and using $\bb<\bt$, (ix) follows. To see (x) and (xi), use (ii), (v) and (vii) to get
\ben
\frac{v^{\prime\prime}(r)}{(v^\prime(r))^2}=\frac{ \bt-1+(\bb-1) r^{p\bt} }{ \bt r^\bt} = \frac{ \bt-1 }{ \bt r^\bt}+ \frac{ \bb-1}{\bt r^{\bb} }.
\een
Since $\bb<\bt$ and $r>1$, the estimates in (xi) hold.
\end{proof}
\vsp

\begin{rem}\label{sec3.15} We now list observations based on Lemma \ref{sec3.10}. These arise from the various cases
described in Theorems \ref{sec2.10} and \ref{sec2.11}.
 
Recall that $k=k_1+1$, $\g=k+1=k_1+2$ and $\s$ is as in Theorems \ref{sec2.10} and \ref{sec2.11}. 
Set $\gs=\g/k.$ We discuss the following three cases.
\ben
&&\mbox{Case (A)\quad $\bt=\gs=2$ and $\bb=2-\ve$, where $0<\ve<1$ and $k=1$.}\\
&&\mbox{Case (B)\quad $\bt=\bb=\gs$ and $k>1$.}\\
&&\mbox{Case (C)\quad $\bt=\gs$ and $\bb=\s/(\s-1)$, where $\s>\g$ and $k\ge 1$.}
\een
\vsp
{\bf Case (A) $k=1$:} Take $\bt=2$ and $\bb=2-\ve$. From Lemma \ref{sec3.10}, $p=\ve/2$ and 
$$v(r)=\int_0^{ r^2 } (1+\tau^{\ve/2})^{-1} \;d\tau. $$

Let $0<\ve<1$, then $1-p=(2-\ve)/2>0$. We apply Lemma \ref{sec3.10} (iii), (iv), (vi), (vii), (viii), (ix) and (xi).
Thus,
\ben
&&(iii)\quad \frac{ \min\left( r^{2-\ve},\;r^2\right) }{2}\le \frac{ r^2}{ 1+ r^\ve}\le v(r)\le \min(r^{2-\ve},\;r^2), \;\forall r \ge 0,\\
&&(iv)\quad \frac{r^{2-\ve} -R^{2-\ve} }{2}\le v(r)-v(R)\le 2\left(r^{2-\ve} -R^{2-\ve} \right), \;\mbox{$\forall r\ge R$ and $\forall R>1$.}
\een

Next,
\ben
(vi)\;  \min\left( r^{2 -\ve }, \; r^2\right)\le rv^{\prime}(r)\le 2 \min\left( r^{2-\ve }, \; r^2\right),\;\;
(viii)\; \frac{v^{\prime}(r)}{r}\le 2\min(1,\;r^{-\ve}).
\een

Finally,
\ben
(ix)\;\;1-\ve\le \frac{rv^{\prime\prime}(r)}{v^{\prime}(r)}\le 1,\;\forall r>0,\;\;\;
(xi)\;\; \frac{1-\ve}{2 r^{2-\ve}}\le \frac{v^{\prime\prime}(r)}{(v^\prime(r))^2}\le \frac{1}{r^{2-\ve} },\;\;\forall r>1.
\een
\vsp
{\bf Case (B) $k>1$:} Set $\bt=\bb=\gs$ and $v(r)=r^{\gs}.$
\vsp
Using that $\g=k+1$ and $k(\gs-1)=1$, we have
\ben
&&(vii)\quad rv^{\prime}(r)=\gs r^{\gs},\quad (viii)\quad \frac{ (v^{\prime}(r))^k }{r}=\left(\gs \right)^k,\quad 
(ix)\quad \frac{r v^{\prime\prime}(r)}{v^{\prime}(r)}=\gs-1=\frac{1}{k},\\
&&(xi)\quad \frac{ v^{\prime\prime}(r)}{( v^{\prime}(r) )^2}=\left( \frac{\gs-1}{\gs} \right)r^{-\gs}=\frac{1}{\g r^{\gs}}.
\een
\vsp
{\bf Case (C) $k\ge 1$:} Set $\bt=\gs$ and $\bb=\s/(\s-1),$ where $\s>\g.$
\vsp
Since $\s>\g$, we have that $\bt>\bb$. Using that $\g=k+1$, we get 
\ben
p=\frac{\bt-\bb}{\bt}=\frac{ \gs- \s/(\s-1) }{\gs}=\frac{\g(\s-1)-k\s}{\g(\s-1)}=\frac{\s -\g}{\g(\s-1)}>0.
\een
Set 
$$v(r)=\int_0^{r^{\gs}} (1+\tau^p)^{-1}\;d\tau.$$
 
We list the observations obtained by applying parts (iii), (iv), (vii), (viii), (ix) and (xi) of Lemma \ref{sec3.10}.

Let $R>1$. Parts (iii) and (iv) read
\ben
&&(iii)\;\; \frac{\min\left( r^{\gs},\;r^{\s/(\s-1)} \right) }{2}\le v(r)\le \min( r^{\gs},\;r^{\s/(\s-1)} ),\;\;\forall r\ge 0,\\
&&(iv) \;\; c_p\left( r^{\s/(\s-1)} -R^{\s/(\s-1)} \right)\le v(r)-v(R)\le (2c_p)\left( r^{\s/(\s-1)} -R^{\s/(\s-1)} \right),\;\forall r\ge R.
\een

Next,
\ben
&&(vii)\;  \frac{\gs \min\left( r^{\s/(\s-1)},\;r^{\gs}\right) }{ 2}\le rv^{\prime}(r)\le \gs\min\left( r^{\s/(\s-1)},\;r^{\gs}\right)
,\nonumber\\
&&(viii)\qquad\frac{ (v^{\prime}(r))^k }{ r}=  \left( \frac{ \gs}{ 1+r^{p\gs} } \right)^k,\quad \mbox{and}\\
&&\left( \frac{\gs}{2}\right)^k\min\left( 1, \;\frac{1}{r^{(\s -\g)/(\s-1)}} \right)\le \frac{ (v^{\prime}(r))^k }{ r}\le(\gs)^k\min\left( 1, \;\frac{1}{ r^{(\s -\g)/(\s-1)} }\right).
\een
The lower bounds in (iii), (vii) and (viii) have been obtained by considering the intervals $0\le r\le 1$ and $r\ge 1$.
\vsp
Finally, since $\s>\g\ge 2$, Lemma \ref{sec3.10} (ix) and (xi) read
\ben
&&(ix)\;\frac{1}{\s}\le \frac{rv^{\prime\prime}(r)}{v^{\prime}(r)}\le \gs-1,\;\;\;(xi)\;\;\frac{ (\gs \s)^{-1}}{ r^{\s/(\s-1)}}\le \frac{v^{\prime\prime}(r)}{(v^\prime(r))^2}\le
 \frac{2}{r^{\s/(\s-1)} },\;\forall r\ge 1. \;\;\;\Box
\een
\end{rem}
\vsp
We make an observation that applies to the various auxiliary functions we make use of in this work.
\vsp
\begin{rem}\label{sec3.19} The sub-solutions and super-solutions in this work involve a $C^1$ function of $t$ and a $C^1$ function $v(r)$, see Remark \ref{sec3.15}. We verify that the expressions for the operator $H$, that arise from the use of these functions, hold in the sense of viscosity at $r=0$. For $r>0$, $v(r)$ is $C^\infty$. See Lemma \ref{sec3.10}.

Let $\kappa(t)\ge 0$ be a $C^1$ function in $t\ge 0$. Set $r=|x|$ and $w(x,t)= \kappa(t) v(r)$, where $v(r)$ is as in (B) and (C) in Remark \ref{sec3.15}. Note that in (A), $v$ is $C^2$. Thus, we discuss 
$$v(r)=\left\{ \begin{array}{ccc} r^{\gs}, & \bb=\bt=\gs,\\  \int_0^{r^{\gs}} (1+\tau^p)^{-1} d\tau,&\quad \bb=\gs,\;\bt=\s(\s-1)^{-1}. \end{array} \right. $$
Here $\gs=\g/k.$ Since $k\ge 1$, we have that $1\le \gs\le 2$. The case of interest is $\gs<2$.
\vsp

Recall (\ref{sec3.4}) in Remark \ref{sec3.2}. Taking $r>0$ and setting $e=x/r$ and $w=\kp(t) v(r)$, we get with a slight rearrangement
\bea\label{sec3.190}
&&H(Dw, D^2w+Z(w) Dw\otimes Dw)+\chi(t) |Dw|^\s-w_t=\chi(t)(\kp(t))^\s|v^\prime(r)|^\s-\kp^\prime(t)v(r) \nonumber\\
&&\quad\qquad\qquad+\frac{(v^{\prime}(r)\kappa(t))^k} {r}   H\left( e, \;I +\left( \frac{rv^{\prime\prime}(r)}{v^{\prime}(r)}-1 \;+\;\kappa(t)(rv^{\prime}(r)) Z(w) \; \right)e\otimes e\;\right).
\eea
\vsp
We now recall parts (viii) and (ix) in Lemma \ref{sec3.10}, See also Cases (B) and (C) in Remark \ref{sec3.15}. Thus, $v(0)=v^\prime(0)=0$ and $rv^{\prime\prime}(r)/v^\prime(r)\rightarrow \gs-1$ and
$(v^\prime(r))^k/r\rightarrow (\gs)^k$ as $r\rightarrow 0$. It is clear that the right hand side of (\ref{sec3.190}) may be extended continuously to $r=0$. Set the limit (as $r\rightarrow 0$) of the right hand side of (\ref{sec3.190}) as
\ben
&&\hat{H}(0)+\chi(t)L(\s),\;\;\;\mbox{where}\;\;\hat{H}(0)=(\gs \kappa(s))^k  H\left( e, \;I +(\gs-2)e\otimes e\;\right),\\
&&\mbox{and $L(\s)=1$, if $\s=0$, and $L(\s)=0,$ if $\s=0$.}
\een
Note that $\hat{H}(0)\ge 0$ since $\gs-2\ge -1.$
Our goal is to show that
\eqRef{sec3.191}
H(Dw, D^2w+Z(w) Dw\otimes Dw)+\chi(t) |Dw|^\s-w_t=\hat{H}(0)+\chi(s) L(\s). 
\ee
holds at points $(0,s)$, i,e, at $r=0$ and $s>0$, in the viscosity sense.

Let $s>0$. Suppose that $\psi$, $C^1$ in $t$ and $C^2$ in $x$, is such that $(w-\psi)(x,t)\le (w-\psi)(o,s)$, for $(x,t)$ near $(o,s)$. Thus,
$$w(x,t)=\kp(t)v(r)\le \langle D\psi(o,s), x\rangle+\psi_t(o,s)(t-s)+o(|x|+|t-s|),$$
as $(x,t)\rightarrow (o,s)$. Since $\kp(t),\;v(r)\ge 0$ and $v^\prime(0)=0$, we have that $\psi_t(o,s)=0$, $D\psi(o,s)=0$. 

Next, using $w(x,t)=\kp(t)v(r)\le \langle D^2\psi(o,s)x, x\rangle/2+o(|x|^2+|t-s|),$
as $(x,t)\rightarrow (o,s)$, and recalling $\bt<2$ and Lemma \ref{sec3.10}(vii), it is clear that $D^2\psi(o,s)$ does not exist and $w$ is a sub-solution of (\ref{sec3.190}). 
\vsp

Now, let $\psi$, $C^1$ in $t$ and $C^2$ in $x$, be such that $(w-\psi)(x,t)\ge (u-\psi)(o,s)$, for $(x,t)$ near $(o,s)$. Thus,
$\kp(t)v(r)\ge \langle D\psi(o,s), x\rangle+\psi_t(o,s)(t-s)+o(|x|+|t-s|),$ as $(x,t)\rightarrow (o,s)$. As argued before,
$D\psi(o,s)=0$ and $\psi_t(o,s)=0$. If $D^2\psi(o,s)$ does not exist then $w$ is a super-solution. 

If $D^2\psi(o,s)$ exists then
\bea\label{sec3.192}
&&H\left(D\psi, D^2\psi+ Z(\psi) D\psi\otimes D\psi \right)(o,s)+\chi(s)|D\psi|^\s(o,s)-\psi_t(o,s) \nonumber\\
&&\qquad\qquad\qquad=H\left(0, D^2\psi \right)(o,s)+\chi(s)L(\s). 
\eea
We now observe that since $\bt=\g/k<2$ and $\g=k+1$, we have $k=k_1+1>1$ and hence, $k_1>0$. Applying Condition B (see(\ref{sec2.2})), $H(0, D^2\psi)(o,s)=0$ and (\ref{sec3.192}) reads
$$H\left(D\psi, D^2\psi+ Z(\psi) D\psi\otimes D\psi \right)(o,s)+\chi(s)|D\psi|^\s(o,s)-\psi_t(o,s)=\chi(s)L(\s).$$
Since $\hat{H}(0)\ge 0$, using (\ref{sec3.191}), we see that $w$ is a super-solution. $\Box$
\end{rem}

\section{\bf Super-solutions}

Our goal in this section is to construct super-solutions whose growth rates, for large $r$, are as stated in Theorem \ref{sec2.10}. Tthe auxiliary functions discussed in Remark \ref{sec3.15} are used to achieve our goal. The construction involves making separate estimates for small $r$ and for large $r$. For small $r$, we employ (\ref{sec3.4}) and, for large $r$, we use (\ref{sec3.5}), see Remark \ref{sec3.2}. 

The section has been divided into two parts: (I) $0\le \s\le \g$ and (II) $\s>\g$. The work in Part I is further divided into two sub-parts (i) $k=1$ and (ii) $k>1$. Part (II) provides a unified work for $k\ge 1$. 
\vsp
The super-solutions we construct are of the kind
\eqRef{sec5.0}
w(x,t)=m+at+b(1+t) v(r),\;\;\mbox{in $\IR^n_T$,}\;\;\;\mbox{where $a\ge 0$ and $0<b<1$,}
\ee
$v$ is $C^1$ in $\IR^n$ and $C^2$ in $\IR^n\setminus\{o\}$ and $-\infty<m<\infty$. We choose $v$ as 
$$\mbox{either}\quad v(r)=\int_0^{r^\bt} \frac{1}{1+\tau^p}\;d\tau\qquad\mbox{or}\qquad v(r)=r^\bt,\;\;\mbox{in $r\ge 0$},$$
for some appropriate $\bt$ and $p$ (or $\bb$), see Lemma \ref{sec3.10}. The scalars $a$ and $b$ are determined later.

In proving that $w$ is a super-solution for appropriate $a$ and $b$, we also calculate the dependence of $a$ on $b$, thus, aiding our calculation of
$\lim_{b\rightarrow 0^+}a$. This is important in showing the claims in Theorem \ref{sec2.10}. 

Throughout this section $\bt=\g/k=\gs$ regardless of the form of $v(r)$, see (\ref{sec5.0}) and Remark \ref{sec3.15}.
The quantity $\bb$, however, depends on $k$ and $\s$, see (\ref{sec5.01}) below.
\vsp
We begin with some preliminary calculations before moving on to Parts I and II. 
Set
\eqRef{sec5.9}
\al=\sup_{[0,T]}|\chi(t)|, \quad \ell=\inf_{\IR}Z,\quad L=\sup_{\IR} Z\quad \mbox{and}\quad \gs=\g/k.
\ee
We assume that $0<\ell\le L<\infty$. We recall that 
\eqRef{sec5.01}
\bt=\gs=\left\{\begin{array}{lcr} 2, & k=1\\ \g/k,&  k>1,\end{array}\right.\quad \mbox{and}\quad \bb=\left\{ \begin{array}{lcr} 2-\ve, & k=1,\;\;0\le \s\le \g,\\ \gs,& k>1,\;\;0\le \s\le \g,\\ \s/(\s-1),& k\ge 1,\;\;\s>\g. \end{array}\right.
\ee
Moreover, we require that
\eqRef{sec5.00}
\mbox{(i)\; if $\s=0$, take $0<\ve <1/8$,}\quad \mbox{and}\quad \mbox{(ii) if $\s>0,$ take $0<\ve<\s/8$.}
\ee
\vsp
Next, we  provide upper bounds for $H$. These will be done for small $r$ and for large $r$ separately. Recall $w$ from (\ref{sec5.0}).
\vsp
\NI{\bf Step 1:}  For small $r$,  we use (\ref{sec3.4}) with $\kappa(t)=1+t$ to obtain that
\bea\label{sec5.2}
&&H(Dw, D^2w+Z(w)Dw\otimes Dw)+\chi(t)|Dw|^{\s}-w_t \nonumber\\
&&\qquad\qquad= \frac{ [ b(1+t) v^{\prime}(r) ]^k}{ r}  H\left(e, I+ \left(  \frac{ r v^{\prime\prime}(r)}{ v^{\prime}(r)  }-1 + b (1+t)(rv^{\prime}(r)) Z(w) \right) e\otimes e \right)\nonumber\\
&&\quad\qquad\qquad\qquad\qquad\qquad +\chi(t) (b(1+t)v^{\prime}(r))^\s-a-b v(r).
\eea
\vsp
For large $r$, we use (\ref{sec3.5})(or (\ref{sec3.7})) to obtain that
\bea\label{sec5.3}
&&H(Dw, D^2w+Z(w)Dw\otimes Dw)+\chi(t)|Dw|^{\s}-w_t \nonumber\\
&&\qquad =b^k[ (1+t) v^{\prime}(r) ]^\g  H\left(e, \frac{I-e\otimes e}{ (1+t) r v^{\prime}(r) }+ \left(  \frac{v^{\prime\prime}(r)}{ (1+t)(v^{\prime}(r))^2  } + bZ(w) \right) e\otimes e \right)
\nonumber\\
&&\quad\qquad\qquad\qquad\qquad+\chi(t)(b (1+t) v^{\prime}(r) )^\s -a - b v(r).
\eea
\vsp
\NI {\bf Step 2:  Bounds for $H$.} 
We employ Remark \ref{sec3.15} and use estimates for $v(r)$(and its derivatives) in (\ref{sec5.2}) and (\ref{sec5.3}) to obtain upper bounds for $H$. Assume $R\ge 1$. A value will be chosen later.

\vsp
{\bf (i) $0\le r\le R$:} Since $Z(w)\le Z(m)\le L$(see (\ref{sec5.9})), define
\eqRef{sec5.4}
M(b,r)=\max_{|e|=1}H\left(e, I+ b(1+T)L \gs r^{\gs}  e\otimes e\right).
\ee
By using the monotonicity in Condition A (see (\ref{sec2.1})(i)) and Condition B (see (\ref{sec2.2})), $M(b, r)$ is non-decreasing in $r$ and $b$, $M(b,r)\ge \max_{|e|=1}H\left(e, I\right)>0$ and
$$M(b,r)\le M(1, R)\le R^{\gs}M(1,1),\;\;\forall R>1.$$
\vsp 
Recall parts (vii) and (ix) of Cases A, B and C in Remark \ref{sec3.15}. It is seen that
\eqRef{sec5.40}
\frac{rv^{\prime\prime}(r)}{v^{\prime}(r)}-1\le \gs-2=\frac{1-k}{k}\le 0\quad\mbox{and}\quad rv^{\prime}\le \gs r^{\gs},\;\;\forall\;k\ge 1.
\ee
We apply the above to (\ref{sec5.2}) and use monotonicity to get
\ben
&&H\left(e, I+ \left(  \frac{ r v^{\prime\prime}(r)}{ v^{\prime}(r)  }-1 + b(1+t) Z(w)r v^{\prime}(r) \right) e\otimes e \right) \nonumber\\
&&\qquad\qquad\qquad\qquad\qquad\qquad\qquad\qquad\le H\left(e, I+ \gs b(1+T)Z(w) r^{\gs}  e\otimes e\right).
\een
Since $Z(w)\le Z(m)\le L$, using (\ref{sec5.4}) and the bound for $M(b,r)$ we obtain that for $0\le r\le R$,
\eqRef{sec5.5}
H\left(e, I+ \left(  \frac{ r v^{\prime\prime}(r)}{ v^{\prime}(r)  } -1+ b(1+t) Z(w)r v^{\prime}(r) \right) e\otimes e \right)\le R^{\gs}M(1,1),\;\;\forall R>1.
\ee
\vsp
Next, we recall the upper bound $(v^{\prime}(r))^k/r\le (\gs)^k$ from part (viii) of Cases A, B and C in Remark \ref{sec3.15}. Thus, (\ref{sec5.9}), (\ref{sec5.2}) and (\ref{sec5.5}) lead to the estimate
\bea\label{sec5.500}
&&H(Dw, D^2w+Z(w)Dw\otimes Dw)+\chi(t)|Dw|^{\s}-w_t \nonumber\\
&&\qquad\qquad\le [b\gs (1+T)]^k  M(1,1) R^{\gs}+\al [b(1+T)v^{\prime}(r)]^\s-a-b v(r).
\eea

\vsp
{\bf (ii) $1\le R\le r$:} We recall parts (vii) and (xi) of Cases A, B and C in Remark \ref{sec3.15}. Part (vii) of the Cases A, B and C show that 
\ben 
rv^{\prime}(r)\ge \left\{ \begin{array}{lcr} r^{2-\ve}, & \quad\;\; k=1,\;\;0\le \s\le \g,\\ \gs r^{\gs}, & \quad\; \; k>1,\;\;0\le \s\le \g,\\ (\gs/2) r^{\s/(\s-1)},& k\ge 1,\;\;\s>\g. \end{array}\right.
\een
In the last inequality, $\s/(\s-1)<\gs$, if $\s>\g$. Thus, using the above and part (xi) of the Cases A, B and C, we obtain
\ben
\max\left( \frac{1}{rv^{\prime}(r)},\;\;\frac{v^{\prime\prime}(r)}{ ( v^{\prime}(r) )^2}\right)\le \left\{ \begin{array}{lcr} 2r^{-(2-\ve)},&\quad \; k=1,\;\;0\le \s\le \g,\\
2 r^{-\gs}, & \quad \; k>1,\;\;0\le \s\le \g,\\
2r^{-\s/(\s-1)}, & k\ge 1,\;\;\s>\g. \end{array}\right.
\een
Thus,
\eqRef{sec5.60}
\max\left( \frac{1}{rv^{\prime}(r)},\;\;\frac{v^{\prime\prime}(r)}{ ( v^{\prime}(r) )^2}\right)\le 2,\;\;\;\mbox{in $r\ge 1$.}
\ee

Noting that both quantities on the left hand side of (\ref{sec5.60}) are non-negative, using Condition A and (\ref{sec5.60}), the term $H$ in (\ref{sec5.3}) yields in $t\ge 0$,
\bea\label{sec5.6}
&&H\left(e, \frac{I-e\otimes e}{ (1+t) r v^{\prime}(r) }+ \left(  \frac{v^{\prime\prime}(r)}{ (1+t)(v^{\prime}(r))^2  } + bZ(w) \right)e\otimes e \right) \nonumber\\
&&\quad\qquad\le H\left(e, \frac{I-e\otimes e}{ r v^{\prime}(r) }+ \left(  \frac{v^{\prime\prime}(r)}{(v^{\prime}(r))^2  } + bZ(w) \right)e\otimes e \right) \nonumber\\
&&\quad\qquad \le H\left(e, 2( I-e\otimes e)+ 2 e\otimes e + bZ(w) I \right)\le H(e, (2+L)I),
\eea
since $I\ge e\otimes e$, $0<b<1$ and $0<Z\le L$.
\vsp
Observing that $w\ge m$, we define
\eqRef{sec5.7}
\bar{M}=\max_{|e|=1} H\left(e, (2+L) I \right).
\ee
Using Conditions A, B and C, $\bar{M}\ge H(e, 2I)=2H(e, I)>0.$
\vsp
Thus, in $r\ge R\ge 1$, by using (\ref{sec5.7}) in (\ref{sec5.6}) we get 
\ben
H\left(e, \frac{I-e\otimes e}{ (1+t) r v^{\prime}(r) }+ \left(  \frac{v^{\prime\prime}(r)}{ (1+t)(v^{\prime}(r))^2  } + Z(w)b \right) e\otimes e \right)\le \bar{M}.
\een
Using (\ref{sec5.9}) and the above upper bound in (\ref{sec5.3}) we get
\bea\label{sec5.8}
&&H(Dw, D^2w+Z(w)Dw\otimes Dw)+\chi(t)|Dw|^{\s}-w_t \nonumber\\
&&\qquad\qquad \le b^k[ (1+T) v^{\prime}(r) ]^\g  \bar{M}+\al [b (1+T) v^{\prime}(r) ]^\s -a - b v(r).
\eea 
\vsp
\NI{\bf Step 3: Additional bounds:}
We record the following bounds that would be useful for what follows. Refer to part (vii) of Cases A, B and C in Remark \ref{sec3.15}. In $r\ge 0$,
\eqRef{sec5.70}
v^{\prime}(r)\le \left\{\begin{array}{lcr} 2\min( r^{1-\ve}, r),& k=1,\;\;0\le \s\le \g,\\ \gs r^{\gs-1},& k>1,\;\;0\le \s\le \g,\\ \gs \min( r^{1/(\s-1)}, r^{{\gs}-1}),& k\ge 1,\;\;\s>\g.
\end{array} \right. \qquad \Box
\ee
\vspp
{\bf Constructions of Super-solutions:} We remind the reader that $k_2=1$, $\g=k+1=k_1+2$ and $\gs=\g/k$ throughout.
\vsp
\NI{\bf Part I: $0\le \s \le \g$.} In what follows we take $R\ge 1$, to be determined later.
\vsp
\NI{\bf Sub-part (i): $k=1.$} Thus, $k_1=0$. Let $\ve>0$ be small. Recall from (\ref{sec5.01}) that
$\g=\gs=2$. We take $p=\ve/2$. Thus, using (\ref{sec5.0}) we get
\eqRef{sec5.90}
w(x,t)=m+at+b(1+t)v(r),\;\;\;\mbox{in $\IR^n_T$},
\ee
where
$$v(r)=\int_0^{r^{2} } (1+\tau^{\ve/2})^{-1}\;d\tau,$$
and $a\ge 0$ and $0<b<1$ are to be determined.

We address the interval $0\le r\le R$. Using (\ref{sec5.70}) and $0\le r\le R$, we get that $v^{\prime}(r)\le 2R$. Employing this in the second term on the right hand side of
(\ref{sec5.500}) we get
\ben
&&H(Dw, D^2w+Z(w)Dw\otimes Dw)+\chi(t)|Dw|^{\s}-w_t \\
&&\qquad\qquad\qquad \le   2b(1+T) M(1,1)R^2  +\al ( 2b(1+T))^\s R^\s -a.
\een

We choose 
\eqRef{sec5.91}
a=\left\{ \begin{array} {lcr} \al+2b(1+T) M(1,1)R^2+bR^{2-\ve}/2, &\;\; \s=0,\\ \al ( 2b(1+T))^\s R^\s+2b(1+T) M(1,1)R^2+bR^{2-\ve}/2, &\;\;0<\s\le \g. \end{array}\right.
\ee
This ensures that $w$ is a super-solution in $0\le r\le R$.
\vsp
Next, we address $r\ge R$. We use the estimate $v^{\prime}(r)\le 2 r^{1-\ve}$(see (\ref{sec5.70})) in the second term of the right hand side of (\ref{sec5.8}) to obtain
\bea\label{sec5.290}
&&H(Dw, D^2w+Z(w)Dw\otimes Dw)+\chi(t)|Dw|^{\s}-w_t \nonumber\\
&&\qquad\qquad \le 4b( (1+T) r^{1-\ve})^2  \bar{M}+\al (2b (1+T) r^{1-\ve} )^\s -a - b v(r).
\eea
We apply the lower bound in part (iv) of Case A in Remark \ref{sec3.15}, that is,
$$v(r)\ge \frac{r^{2-\ve} -R^{2-\ve} }{2},\;\;\;\forall r\ge R\ge 1.$$

Thus, we obtain from (\ref{sec5.290}) that
\bea\label{sec5.930}
&&H(Dw, D^2w+Z(w)Dw\otimes Dw)+\chi(t)|Dw|^{\s}-w_t \nonumber\\
&&\quad \le 4b(1+T)^2  \bar{M}r^{2-2\ve}+\al[2 b (1+T)]^\s r^{\s(1-\ve)} -a - b\left( \frac{r^{2-\ve} -R^{2-\ve}}{2} \right) \nonumber\\
&&\quad \le 4b(1+T)^2 \bar{M}  r^{2-2\ve} +\al[2 b (1+T)]^\s r^{\s(1-\ve)} -\hat{a}- \frac{b r^{2-\ve}}{2},
\eea
where in the last inequality we have used the expression for $\hat{a}=a-bR^{2-\ve}/2$, see (\ref{sec5.91}).
\vsp
\NI{\bf (a): $\s=0.$} Using (\ref{sec5.91}) and that $\hat{a}\ge 0$, the right hand side in (\ref{sec5.930}) yields
\ben
4b(1+T)^2 \bar{M}r^{2-2\ve}+\al-a- \left( \frac{ b}{2} \right) r^{2-\ve} \le br^{2-2\ve} \left( 4(1+T)^2\bar{M}- \frac{r^{\ve}}{2} \right).
\een
Choose $R$ such that $R^\ve=\max\left(1,\;8(1+T)^2\bar{M} \right)$. Clearly, $w$ is a super-solution in $\IR^n_T$. 
\vsp
We record that the above choice for $R$ and (\ref{sec5.91}) yield that 
\eqRef{sec5.9300}
\lim_{b\rightarrow 0^+} a=\al, \qquad\mbox{for $\s=0$.}
\ee
\vsp
\NI{\bf (b): $0<\s\le 2$.} Note that $\g=2$ and $\hat{a}\ge 0$. The right hand side of (\ref{sec5.930}) yields
\bea\label{sec5.9301}
&&4b(1+T)^2\bar{M}r^{2-2\ve}+\al[2 b (1+T)]^\s r^{\s(1-\ve)}- \frac{b r^{2-\ve}}{2}-\hat{a}  \nonumber \\
&&\qquad\qquad\le 4b(1+T)^2 \bar{M} r^{2-2\ve} + \al[2b(1+T)]^\s r^{\s(1-\ve)} - \frac{b r^{2-\ve} }{2}   \nonumber\\
&&\qquad\qquad\le b r^{2-2\ve}\left( 4(1+T)^2 \bar{M}+ \frac{ \al[2(1+T)]^\s b^{\s-1} }{ r^{(2-\s)(1-\ve)} } - \frac{R^\ve}{2} \right),
\eea
in $r\ge R\ge 1$. Also, $(2-\s)(1-\ve)\ge 0$.
Set 
\eqRef{sec5.9302}
P=4(1+T)^2 \bar{M}\quad\mbox{and}\quad Q=\al[2(1+T)]^\s.
\ee

Select
\eqRef{sec5.9303}
R=\left\{ \begin{array} {ccr} \max\left\{  ( 2(1+P) )^{1/\ve},\;\; (2 Q b^{\s-1})^{1/(2-\s)(1-\ve)} \right\},& 0<\s<1, \\
\max\left\{1,\;\;( 2P+2Q)^{1/\ve} \right\},& 1\le \s\le 2.
 \end{array} \right.
\ee
For $0<\s<1$, we have set $r=R$ in the second term of (\ref{sec5.9301}) and chosen $R$, and 
for $1\le \s\le 2$, we have taken $r=b=1$ in the second term of (\ref{sec5.9301}).

Using (\ref{sec5.9302}) and (\ref{sec5.9303}) in (\ref{sec5.9301}) and recalling (\ref{sec5.930}), $w$ is a super-solution in $\IR^n_T$.
\vsp
We recall the expression for $a$ in (\ref{sec5.91}) and claim that $a\rightarrow 0$ as $b\rightarrow 0$. This is clear for $1\le \s\le 2$ because of the choice in
(\ref{sec5.9303}). 
The case of interest is $0<\s<1$ since $R\rightarrow \infty$ as $b\rightarrow 0$. It suffices to show that $bR^2\rightarrow 0$ as $b\rightarrow 0$ as
this would imply the same of $bR$ and $bR^{2-\ve}$. Taking $b$ small in (\ref{sec5.9303}), one can write 
$$R=K b^{(\s-1)/(2-\s)(1-\ve)}\qquad \mbox{and}\qquad bR^2=K^2 b^{1+2(\s-1)/(2-\s)(1-\ve)},$$
for an appropriate $K$ that is independent of $b$. A simple calculation shows that
\ben
1+ \frac{ 2(\s-1)}{(2-\s)(1-\ve)}=\frac{\s(1+\ve)-2\ve}{(2-\s)(1-\ve)}.
\een
From (\ref{sec5.00}), $0<\ve<\s/8$ and, hence, $\s(1+\ve)-2\ve>0$. Recalling (\ref{sec5.90}), (\ref{sec5.91}) and (\ref{sec5.9300}) we obtain that, for any small $\ve>0$,
\eqRef{sec5.94}
\lim_{b\rightarrow 0} a=\left\{ \begin{array}{lcr} \al, & \quad \s=0,\\  0,&\quad 0<\s\le 2. \end{array}  \right.
\ee
\vsp
\NI{\bf Sub-part (ii): $k>1$.} We set
\eqRef{sec5.95}
w(x,t)=m+at+b(1+t) r^{\gs},\;\;\mbox{in $\IR^n_T$.}
\ee
Recall that $\gs=\g/k=1+1/k$, $k=k_1+1$ and $\g=k+1=k_1+2$.
\vsp
Consider $0\le r\le R$, where $R>1$ is to be determined. We recall (\ref{sec5.500}) and use $v^{\prime}(r)=\gs r^{\gs-1}$ to obtain
\ben
&&H(Dw, D^2w+Z(w)Dw\otimes Dw)+\chi(t)|Dw|^\s-w_t \\
&&\qquad\qquad\qquad \le [ b\gs (1+T) ]^k M(1,1)R^{\gs}+\al [b\gs (1+T) R^{\gs-1}]^\s-a-bv(r).
\een

Noting that $\gs-1=1/k$, we choose,
\eqRef{sec5.96}
a=b^k[\gs (1+T) ]^k M(1,1)R^{\gs}+\al b^\s[\gs (1+T)]^\s R^{\s/k}.
\ee
Thus, $w$ is a super-solution in $\overline{B_R(o)}\times (0,T).$
\vsp
Now consider $r\ge R$. Using (\ref{sec5.8}) and $v^{\prime}(r)=\gs r^{\gs-1}$, we get
\ben
&&H(Dw, D^2w+Z(w)Dw\otimes Dw)+\chi(t)|Dw|^\s-w_t \nonumber\\
&&\qquad\qquad\qquad \le b^k [\gs (1+T) r^{\gs-1} ]^\g  \bar{M}+\al [b\gs (1+T) r^{\gs-1}]^\s-a- b r^{\gs}.
\een
Since $\gs-1=1/k$, clearly, $\g(\gs-1)=\gs$. Setting $E=\gs (1+T)$, the above reads 
\bea\label{sec5.97}
&&H(Dw, D^2w+Z(w)Dw\otimes Dw)+\chi(t) |Dw|^\s-w_t \nonumber\\
&&\qquad\qquad\qquad\qquad\qquad \le \left(E^\g \bar{M}\right) b^kr^{\gs}+ \left(\al E^\s \right) b^\s r^{\s/k}-a-br^{\gs}.
\eea
\vsp
We analyze separately: (1) $\s=0$, (2) $1<\s\le \g$, and (3) $0<\s\le 1$.
\vsp
{\bf (1) $\s=0$:} Setting $R=1$ and recalling (\ref{sec5.96}), the right hand side of (\ref{sec5.97}) yields
\ben
(E^\g \bar{M}) b^kr^{\gs}+\al-a-br^{\gs}\le b r^{\gs} \left( E^\g \bar{M} b^{k-1}-1  \right)\le 0,
\een
if we choose $0<b\le \min\left(1,\; \left( E^\g\bar{M}\right)^{-1/(k-1)} \right).$
Thus, $w$ is a super-solution in $\IR^n_T$ and 
\eqRef{sec5.970}
\lim_{b\rightarrow 0} a=\al.
\ee
\vsp
{\bf (2) $1< \s\le \g$:} The right hand side of (\ref{sec5.97}) is bounded above, in $r\ge R$, by
\eqRef{sec5.971}
(E^\g \bar{M}) b^kr^{\gs}+(\al E^\s) b^\s r^{\s/k}-a-br^{\gs}\le b r^{\gs}\left[ (E^k \bar{M}) b^{k-1}+ \frac{ (\al E^\s) b^{\s-1} }{ R^{(\g-\s)/k} } -1\right].
\ee

Setting $R=1$ in the second term of the right hand side of (\ref{sec5.971}), we get
\ben
(E^\g \bar{M}) b^kr^{\gs}+(\al E^\s) b^\s r^{\s/k}-br^{\gs}\le b r^{\gs}\left[ (E^k \bar{M}) b^{k-1}+(\al E^\s) b^{\s-1} -1\right].
\een
Choosing $0<b<1$, small enough, we get that $w$ is super-solution in $\IR^n_T$. Moreover, using (\ref{sec5.96}) $\lim_{b\rightarrow 0} a=0$.
\vsp 
{\bf (3) $0<\s\le 1$:} We recall (\ref{sec5.971}) i.e.,
\ben
(E^\g \bar{M}) b^kr^{\gs}+(\al E^\s) b^\s r^{\s/k}-a-br^{\gs}\le b r^{\gs}\left[ (E^\g \bar{M}) b^{k-1}+ \frac{ (\al E^\s) b^{\s-1} }{ R^{(\g-\s)/k} } -1\right].
\een
We choose
$$b<\min\left[1,\; \left( \frac{1}{4 E^\g \bar{M} } \right)^{1(k-1)} \right]\quad\mbox{and}\quad R=\max\left[ 1, \;    \left\{ (4\al E^\s) b^{\s-1} \right\}^{k/(\g-\s)}    \right].$$
It is clear that $w$ is a super-solution in $\IR^n_T$.

Our next task is to show that $\lim_{b\rightarrow 0} a=0$. Recalling (\ref{sec5.96}) and comparing the terms $b^k R^{\gs}$ and $(bR^{1/k})^\s$, we see that it is enough to show that $b^k R^{\gs}\rightarrow 0$ as $b\rightarrow 0$. This is clear if $\s=1$. Assuming that $\s<1$ and using the choice for $R$, we
see that (use $\gs=\g/k$)
$$b^k R^{\gs}=K b^k  \left[ b^{k(\s-1)/(\g- \s)}\right]^{\gs}=K b^{k+\g(\s-1)/(\g-\s)},$$
for some $K$ independent of $b$. Using that $\g=k+1=k_1+2$, we calculate
$$k+\frac{ \g(\s-1) }{ \g-\s}=\frac{k_1\g+\s}{\g-\s}>0.$$
The claim holds. 

Summarizing from Sub-Parts (i) (see (\ref{sec5.94})) and (ii) (see (1), (2) and (3)), we get 
\eqRef{sec5.98}
\lim_{b\rightarrow 0} a=\left\{ \begin{array}{lcr} \al, &\quad \s=0,\\ 0, & \quad 0<\s\le \g.  \end{array}  \right.
\ee
\vspp
{\bf Part II $\s> \g,\;k\ge 1$:}  We set
\bea\label{sec5.99}
&&w(x,t)=m+at+b(1+t) v(r),\;\;\mbox{where}\nonumber\\
&&\quad v(r)=\int_0^{ r^{\gs} }  \frac{1}{ 1+\tau^p} d\tau\;\;\;\mbox{and}\;\;p= \frac{ \s-\g }{ \g(\s-1) }.
\eea
We recall estimates stated in Case C of Remark \ref{sec3.15}.
 \vsp
Take $R\ge 1$ and consider $0\le r\le R$. We employ (\ref{sec5.500}) i.e.,
\ben
&&H(Dw, D^2w+Z(w)Dw\otimes Dw)+\chi(t)|Dw|^{\s}-w_t \nonumber\\
&&\qquad\qquad\qquad\le [b\gs (1+T)]^k M(1,1)R^{\gs}+\al [b(1+T)v^{\prime}(r)]^\s-a-b v(r).
\een
Noting that $\s-1\ge \g-1=k$, using (\ref{sec5.70}) ($v^{\prime}(r)\le \gs r^{1/(\s-1)}$) and setting $E=\gs (1+T)$, we get from above that
\ben
&&H(Dw, D^2w+Z(w)Dw\otimes Dw)+\chi(t)|Dw|^{\s}-w_t\\
&&\qquad\qquad\qquad \le (b E)^k  M(1,1)R^{\gs}+\al (b E)R^{\s/(\s-1)}-a.
\een

Set $c_p=[2(1-p)]^{-1}$ and select
\eqRef{sec5.100}
a=(bE)^k  M(1,1)R^{\gs}+\al (bE)^\s R^{\s/(\s-1)}+ c_p b R^{ \s/(\s-1)}.
\ee
Thus, $w$ is a super-solution in $B_R(o)\times (0,T)$.
\vsp
In $r\ge R$, we use (\ref{sec5.8}) i.e.,
\bea\label{sec5.101}
&&H(Dw, D^2w+Z(w)Dw\otimes Dw)+\chi(t)|Dw|^{\s}-w_t \nonumber\\
&&\qquad\qquad\qquad \le b^k[ (1+T) v^{\prime}(r) ]^\g  \bar{M}+\al [b (1+T) v^{\prime}(r) ]^\s -a - b v(r).
\eea 
From part (iv) of Case C in Remark \ref{sec3.15}, we have
\ben
v(r)\ge c_p \left( r^{ \s/(\s-1) } - R^{ \s/(\s-1) } \right) ,\;\;\forall r\ge R.
\een
Using (\ref{sec5.70}) ($v^{\prime}(r)\le \gs r^{1/(\s-1)}$), the lower bound for $v(r)$ stated above, $E=\gs (1+T)$ and (\ref{sec5.100}) in the right hand side of (\ref{sec5.101}), we get
\bea\label{sec5.1010}
&&b^k [ (1+T) v^{\prime}(r) ]^\g  \bar{M}+\al [b (1+T) v^{\prime}(r) ]^\s -a - b v(r) \nonumber\\
&&\quad\qquad \le b^k(E^{\g} \bar{M}) r^{\g/(\s-1) }  + b^\s(\al E^\s) r^{ \s/(\s-1) } -a-c_p b \left( r^{\s/(\s-1)} -R^{ \s/(\s-1) } \right) \nonumber\\
&&\quad\qquad \le b^k (E^{\g} \bar{M}) r^{\g/(\s-1) }  + b^\s(\al E^\s) r^{ \s/(\s-1) } -c_p b r^{\s/(\s-1)}  \nonumber\\
&&\quad\qquad \le b r^{\s/(\s-1)}\left[ \frac{  b^{k-1} (E^{\g} \bar{M}) }{ R^{(\s-\g)/(\s-1) } } + b^{\s-1}(\al E^\s) -c_p \right],
\eea
where we have used $1<\g<\s$ and $r\ge R$.
\vsp
For $k>1$, take $R=1$ and $b>0$ small enough (depending on $\s,\;\al,\;E$ and $\bar{M}$) so that (\ref{sec5.1010}) is negative. 
If $k=1$ we take
\eqRef{sec5.103}
R=\max\left[ 1,\;  \left(  \frac{ 4 E^{\g} \bar{M} }{ c_p }    \right)^{(\s-1)/(\s-\g)} \right] \quad \mbox{and}\quad b \le \min \left[1,\; \left( \frac{ c_p }{ 4\al E^\s } \right)^{ 1/(\s-1) } \right].
\ee
With these selections, the right hand side of (\ref{sec5.1010}) is negative. Thus, (\ref{sec5.101}) implies that $w$ is super-solution in $\IR^n_T$. Recalling (\ref{sec5.100}), we see that
\eqRef{sec5.200}
\lim_{b\rightarrow 0} a=0.
\ee

\section{\bf Sub-solutions}

In this section, we construct sub-solutions. We place no restrictions on the growth rate if $0\le \s<\g$. This includes also the case when $\sup_{[0,T]}|\chi(t)|$ is small enough.
However, in general, a lower bound in the case $\s\ge \g$ is needed for our work. We remark that the auxiliary functions employed are closely related to the functions used for super-solutions. 

We achieve our goal by utilizing the expressions in Remark \ref{sec3.2}, in particular, the versions in (\ref{sec3.8}) and (\ref{sec3.80}).  
Thus, setting $w(x,t)=v(r)-\kappa(t)$ and assuming that $v^{\prime}(r)\le 0,$ we get that
\bea\label{sec6.18}
&&H(Dw, D^2w+Z(w) Dw\otimes Dw)+\chi(t) |Dw|^\s-w_t \nonumber\\
&&\qquad\qquad\qquad=\frac{|v^\prime(r)|^k } {r }H\left(e, \left(r|v^{\prime}(r)| Z(w)+1- \frac{rv^{\prime\prime}(r)}{v^{\prime}(r)} \right) e\otimes e-I\right)\nonumber\\
&&\qquad\qquad\qquad\qquad\qquad+ \chi(t)|v^{\prime}(r)|^\s+\kappa^{\prime}(t).
\eea
\vsp
Next, we recall Condition C (see (\ref{sec2.5})), (\ref{sec2.51}) and (\ref{sec3.0}) and set
\eqRef{sec6.180}
N=\min_{|e|=1}H(e,-I),\quad K_0=\frac{\Lambda(\lam_0)}{\lam_0}\quad \mbox{and}\quad \ell=\inf_{s}Z(s).
\ee
Set $\cH(\lam)=\min_{|e|=1}H(e, e\otimes e-\lam^{-1}I)$ and $\cH=\min_{|e|=1}H(e, e\otimes e).$ We record that
\eqRef{sec6.181}
0<\ell<\infty,\quad N<0,\quad \cH(\lam)\ge K_0>0,\;\forall \lam\ge \lam_0,\quad\mbox{and}\quad \lim_{\lam\rightarrow \infty}\cH(\lam)=\cH.
\ee
\vsp
{\bf An auxiliary function and preliminary calculations.}
\vsp
Fix $R>1$. Let $p\ge 1$ and $E\ge 0$, to be determined later. In $0\le r<R$, set 
\eqRef{sec6.19}
\om=r/R,\;\;\;v(\om)=E\int_0^{\om^2} (\tau^p-1)^{-1}d\tau=E R^{2(p-1)} \int_0^{r^2} (\tau^p-R^{2p})^{-1} d\tau.
\ee
Hence, $v$ is defined in $0\le \om<1$. We will often write $v(\om)$ as $v(r)$.

Clearly,
$$v\le 0,\quad v(0)=0,\quad v^{\prime}(r)\le 0\;\;\mbox{and}\;\;\;v(r)\rightarrow -\infty\;\;\mbox{as $r\rightarrow R$.}$$
Set 
\eqRef{sec6.200}
L=L(\om)=\frac{2E}{ 1-\om^{2p} },\;\;\forall 0\le \om<1.
\ee

Differentiating $v(r)$ in (\ref{sec6.19}) and using (\ref{sec6.200}), we get
\bea\label{sec6.20}
&&(i)\;v^{\prime}(r)=-\frac{L(\om)\om}{R}=-\frac{L(\om) r}{R^2} ,\quad 
(ii)\;v^{\prime\prime}(r)=-\frac{L(\om) }{R^2} \left(  \frac{ 1+(2p-1) \om^{2p} } {1-\om^{2p}   } \right),\nonumber\\
&&\quad(iii)\;\frac{ r v^{\prime\prime}(r) }{ v^{\prime}(r) }= \frac{1 +(2p-1) \om^{2p} }{ 1-\om^{2p} },\quad\mbox{and}
\quad(iv)\;\frac{ r v^{\prime\prime}(r) }{ v^{\prime}(r) }-1=\frac{2p \om^{2p} }{ 1 - \om^{2p} }.
\eea
\vsp
Using $k=k_1+1$, $\g=k_1+2$ and (\ref{sec6.20})(i), we get
\eqRef{sec6.21}
\frac{|v^{\prime}(r)|^k}{r}=\frac{ L(\om)^k \om^k }{ R^{k}(\om R) }=\frac{L(\om)^k\om^{k_1} }{R^\g}.
\ee

Next, recalling (\ref{sec6.18}), (\ref{sec6.181}), (\ref{sec6.200}) and (\ref{sec6.20})(i) and (iv), we see that
\bea\label{sec6.22}
&&1- \frac{rv^{\prime\prime}(r)}{v^{\prime}(r)}+rZ(\cdot) |v^{\prime}(r)|=Z(\cdot)L(\om)\om^2-  \frac{2p \om^{2p} }{ 1-\om^{2p} } =\frac{  2EZ(\cdot) \om^2  }{1-\om^{2p} }-\frac{2p \om^{2p} }{ 1-\om^{2p} }   \nonumber\\    
&&\quad\qquad\qquad\qquad\qquad\qquad\ge 2\om^2 \left( \frac{\ell E-p \om^{2(p-1)} }{ 1-\om^{2p} }\right).
\eea
\vsp
\NI{\bf Sub-solutions.}
\vsp
We provide separate treatments for $0\le \s\le \g$ and $\s\ge \g$. The case $\s=\g$ will be addressed in both situations.
\vsp
\NI{\bf Case I: $0\le \s\le\g$.}  Let $\mu\in (-\infty,\infty)$ and recall (\ref{sec6.19}). Set in $0\le r<R$, $\om=r/R$,
\eqRef{sec6.23}
\bar{w}(x,t)=\mu+v(r)-F t, \;\;\mbox{where}\;\;\;v(r)=E\int_0^{\om^2} (\tau^p-1)^{-1} d\tau,
\ee
where $E$, $F$ and $p\ge 2$ are to be determined. Of importance is the limit $\lim_{R\rightarrow\infty}F(R).$
\vsp
Employing (\ref{sec6.18}), (\ref{sec6.200}), (\ref{sec6.21}) and (\ref{sec6.22}), we see that
\bea\label{sec6.24}
&&H(D\wb, D^2\wb+Z(\wb) D\wb\otimes D\wb)+\chi(t) |D\wb|^\s-\wb_t\nonumber\\
&&\qquad\ge  \frac{  L(\om)^k  \om^{k_1} }{ R^{\g} }  H\left( e, 2\om^2\left( \frac{ \ell E-p \om^{2(p-1)} }{1-\om^{2p} } \right)e\otimes e-I  \right)-\al  \left(\frac{ L(\om)\om  }{R } \right)^\s+F.
\eea
\vsp
The sub-solution we construct will depend on $p$ and $R$. Select
\eqRef{sec6.241}
E=\frac{p(p+1) }{\ell}\quad\mbox{and}\quad L(\om)=\frac{2p(p+1)}{\ell(1-\om^{2p} )}.
\ee
As $0\le \om<1$ and $p\ge 2$, we get that $2\om^2\left(\ell E -p \om^{2(p-1)} \right)\ge 2p^2\om^2 $. Set
\eqRef{sec6.250}
J_p(\om)=\frac{ 2 p^2 \om^{2} }{1-\om^{2p}}=\left( \frac{p}{p+1}\right)\ell L(\om)\om^2,
\ee
where we have used (\ref{sec6.200}). Recalling (\ref{sec6.24})  we see that 
\bea\label{sec6.25}
&&H(D\wb, D^2\wb+Z(\wb) D\wb\otimes D\wb)+\chi(t) |D\wb|^\s-\wb_t  \nonumber\\
&&\qquad\qquad\qquad\ge  \left( \frac{ L(\om)^k \om^{k_1} }{R^{\g} } \right) H\left( e,  J_p \left(\om  \right)e\otimes e-I  \right)-  \al \left( \frac{ L(\om)\om }{ R } \right)^\s +F.
\eea
\vsp
We fix $1/\sqrt{2}\le \om_0<1$ and consider separately: (i) $0\le \om\le \om_0$, and (ii) $\om_0\le \om <1$. 
\vsp
\NI{\bf (i) $0\le \om\le \om_0 $:}  Recall (\ref{sec6.180}), (\ref{sec6.181}) and (\ref{sec6.250}). 
We bound 
$$H\left( e,  J_p (\om) e\otimes e-I  \right)\ge H(e, -I)\ge -|N|.$$ 
Using the above in (\ref{sec6.25}) we get that 
\bea\label{sec6.251}
&&H(D\wb, D^2\wb+Z(\wb) D\wb\otimes D\wb)+\chi(t) |D\wb|^\s-\wb_t  \nonumber\\
&&\qquad\qquad\qquad\qquad\qquad\qquad\qquad\ge F-\left(\frac{ L( \om )^k |N| \om^{k_1}  }{ R^{\g} }  +   \al \left(\frac{  L( \om) \om }{R } \right)^\s \right).
\eea

From (\ref{sec6.241}), $L(\om)$ is increasing in $\om$. Since $0\le \om\le \om_0$, we choose
\eqRef{sec6.2510}
F= \frac{ L(\om_0)^k   |N| \om_0^{k_1}}{ R^{\g} } +\al \left( \frac{L(\om_0)\om_0 }{  R  }\right)^\s.
\ee
Thus (\ref{sec6.251}) implies that $\wb$ is a sub-solution in $B_{\om_0R}(o)\times (0,T)$.
\vsp
\NI{\bf (ii) $\om_0 \le \om<1$:} The work will lead to a determination of $p$.

We estimate $J_p(\om)$(recall (\ref{sec6.250})). Since $J_p$ is increasing in $\om$, we get that
\eqRef{sec6.252}
J_p(\om)\ge J_p(\om_0)= \frac{ 2 p^2\om_0^2  }{ 1-\om_0^{2p}}\ge p^2 ,
\ee
since $\om_0^2\ge 1/2$. We note also that $J_p(\om_0)\rightarrow \infty$ if $p\rightarrow \infty$.
\vsp
Using the monotonicity and the homogeneity of $H$(Conditions A and B), $\om_0 \le \om<1$ and (\ref{sec6.252}), we have that
\bea\label{sec6.26}
\min_{|e|=1}H\left( e,  J_p \left( \om\right) e\otimes e-I  \right) &\ge& J_p \left( \om \right)\min_{|e|=1}H\left( e,  e\otimes e- \frac{ I }{J_p\left(\om_0) \right) } \right)  \nonumber\\
&\ge& J_p\left(\om\right) \cH(p^2)\ge K_0 J_p(\om_0)\ge K_0 p^2>0.
\eea
Here we have used (\ref{sec6.181}) and chosen $p\ge p_0$, where $p_0\ge 2$ is large enough.
\vsp
From here on we take $p\ge p_0$ such that (\ref{sec6.26}) holds (see (\ref{sec6.252})). 
Next, using (\ref{sec6.26}) in (\ref{sec6.25}), we obtain
\bea\label{sec6.27}
&&H(D\wb, D^2\wb+Z(\wb) D\wb\otimes D\wb)+\chi(t) |D\wb|^\s-\wb_t \nonumber\\
&&\qquad\qquad\qquad\ge \frac{ L(\om)^k \om^{k_1}J_p( \om ) \cH(p^2)  }{ R^{\g} }-\al \left( \frac{L(\om)\om}{ R }\right)^\s +F \nonumber\\
&&\qquad\qquad\qquad=\ell \cH(p^2) \left( \frac{ p}{p+1} \right)  \left( \frac{L(\om)\om}{R} \right)^\g
-\al \left( \frac{ L(\om) \om}{R } \right)^\s+F.
\eea
In the last inequality, we have used (\ref{sec6.250}) and $\g=k+1=k_1+2$.
\vsp
We factor $(\om L(\om)/R)^\s$ from (\ref{sec6.27}) and use that $\om_0\le \om<1$, to obtain that 


\bea\label{sec6.30}
&&H(D\wb, D^2\wb+Z(\wb) D\wb\otimes D\wb)+\chi(t) |D\wb|^\s-\wb_t \nonumber \\
&&\qquad\qquad\qquad\qquad\qquad\ge \left( \frac{\om L}{R} \right)^\s \left[ \ell \cH(p^2) \left( \frac{ p}{p+1} \right)  \left( \frac{L(\om_0)\om_0 }{R} \right)^{\g-\s}   - \al  \right] +F.
\eea
\vsp
We address $0\le \s\le \g$. We make comments about $\s=\g$ in Sub-Case (c).
\vsp
\NI{\bf Sub-Case (a) $0\le \s<\g$:} As noted earlier, $\wb$ is a sub-solution in $B_{\om_0R}(o)\times (0,T)$.

We refer to (\ref{sec6.30}) and select $R$ such that
\eqRef{sec6.31}
\frac{L(\om_0)}{R}=\frac{1}{\om_0}\left[  \left( \frac{\al}{\ell \cH(p^2)} \right) \left( \frac{1+ p}{p }\right) \right]^{1/(\g-\s)}.
\ee
With this choice, $\wb$ is a sub-solution in $B_R(o)\times (0,T)$.

Using (\ref{sec6.241}) and (\ref{sec6.31}), we get that for some $K_1=K_1(\al, \g,\ell, \om_0, K_0)>0$,
$$R= K_1 \left( p^{\g-\s+1} (p+1)^{k-\s} \right)^{1/(\g-\s)}=O(p^2), \;\;\mbox{as $p\rightarrow \infty$},$$
where we have used $\g=k+1$. Thus, $R\rightarrow \infty$ if and only if 
$p\rightarrow \infty$.
\vsp
We now calculate $\lim_{R\rightarrow \infty} F$.
From (\ref{sec6.241}) and (\ref{sec6.2510}), we write $F$ as the sum of two terms $X$ and $Y$ as follows:
\ben
&&F= \frac{ |N| L(\om_0)^k \om_0^{k_1} }{R^{\g} }+\al \left( \frac{L(\om_0)\om_0  }{ R }\right)^\s=X+Y. 
\een
We use (\ref{sec6.31}), $\g=k+1$ and $k=k_1+1$ to observe that
\ben
\lim_{p\rightarrow \infty}X=\lim_{R\rightarrow \infty}X= \frac{|N| \om_0^{k_1}}{R} \left( \frac{L(\om_0)}{R}\right)^k=0.
\een 

Next, using (\ref{sec6.31}), we get
\ben
Y=\al \left( \frac{L(\om_0)\om_0}{R} \right)^\s=\left[   \frac{ \al^\g}{ ( \ell \cH(p^2) )^\s} \left( \frac{p+1}{p} \right)^\s \right]^{1/(\g-\s)}.
\een
Referring to (\ref{sec6.181}), we see that
\eqRef{sec6.32}
\lim_{R\rightarrow \infty} F=\lim_{R\rightarrow \infty}Y=\lim_{p\rightarrow \infty} Y= \left( \frac{ \al^\g}{ (\ell \cH)^\s } \right)^{1/(\g-\s)},     \quad 0\le \s<\g.
\ee
From (\ref{sec6.32}),
\eqRef{sec6.33}
\lim_{R\rightarrow \infty}F=\left\{ \begin{array}{lcr} \qquad \al,&\s=0,\\ ( \al^\g / (\ell \cH)^\s )^{1/(\g-\s)},&    0<\s<\g.  \end{array}\right.
\ee
\vsp
\NI{\bf Sub-Case (b) $\chi\ge 0:$} An inspection of (\ref{sec6.25}), (\ref{sec6.2510}) and (\ref{sec6.30})($-\al$ is replaced by $+\al$) shows that $F=X$. Thus, $\wb$ is a sub-solution in $B_R(o)\times (0,T)$ for any $\s\ge 0$ and any $R>0$, as there are no restrictions on $R$. Clearly, $\lim_{R\rightarrow \infty} F=0.$
\vsp
\NI{\bf Sub-Case (c) $\s=\g$:} An inspection of (\ref{sec6.30}) shows that if 
$$\al< \ell \cH=(\inf_s Z(s))(\min_{|e|=1}H(e, e\otimes e),$$
then by selecting $p$, large enough, the right hand side of (\ref{sec6.27}) may be written as
$$ \left( \frac{\om L}{R} \right)^\g \left[ \ell \cH(p^2) \left( \frac{ p}{p+1} \right)   - \al  \right] +F\ge 0.$$
For the chosen $p$, $\wb$ is a sub-solution in $B_R(o)\times (0,T)$ for any $R>0$. Moreover, $R$ is independent of $p$ and $F(R)\rightarrow 0$ as $R\rightarrow \infty$. However, if $\al$ exceeds the above value then it is not clear if this conclusion holds. See Case II below.   
\vsp
\NI{\bf Case II $\g\le \s<\infty$:} We assume a lower bound for $u$ and adapt the work in Section 5.  
See also the bounds on $H$ which appear in the beginning of Section 5.

Recall that $k_2=1$, $\g=k+1=k_1+2$ and $\gs=\g/k$. We divide the work into two sub-cases.
\vsp
\NI{\bf Sub-Case (i) $\s=\g:$} We assume that $\al\ge \ell \cH$ and refer to Sub-Parts (i) and (ii) of Part I in Section 5. 
\vsp
\NI{\bf (i1) $k=1$:} Here $\g=\gs=2$. We assume that for any $\ve>0$, small, $\sup_{|x|\ge r} (-u(x,t))\le o(|r|^{2-\ve})$ as $r\rightarrow \infty$. We take
$$\wb(x,t)=m-at-b(1+t)v(r),\quad\mbox{where}\quad v(r)=\int_0^{r^2} (1+\tau^{\ve/2} )^{-1} d\tau.$$
\vsp
\NI{\bf (i2) $k>1$:} Thus, $\g=k+1>2$. We assume that $\sup_{|x|\ge r} (-u(x,t))\le o(|r|^{\gs})$ as $r\rightarrow \infty$. We take
$$\wb=m-at-b(1+t)r^{\gs}.$$
\vsp
\NI{\bf Sub-Case (ii) $\s>\g:$} We allow $k\ge 1$ and refer to Part II of Section 5. We assume that $\sup_{|x|\ge r} (-u(x,t))\le o(|r|^{ \s/(\s-1) } )$ as $r\rightarrow \infty$.
We take
$$\wb=m-at-b(1+t) v(r),\;\;\mbox{where}\;\;v(r)=\int_0^{r^{\gs} } (1+\tau^p)^{-1} d\tau\;\; \mbox{and}\;\;p=\frac{\s-\g}{\g(\s-1)}.$$
\vsp
Since $\wb_r\le 0$ for all the cases described above, we recall the two versions in (\ref{sec3.8}), i.e., for $R>0$, to be determined,
\bea\label{sec6.34}
&& H(D\wb, D^2\wb+Z(\wb) D\wb\otimes D\wb)     \\
&&\qquad=\frac{|\wb_r|^k}{ r } H\left( e, \;\left(\; r|\wb_r|Z(\wb)\;+1   - \frac{r\wb_{rr}}{\wb_r}\right)e\otimes e -I \right),\;\;  \forall 0\le r\le R, \nonumber\\
&&\qquad=|\wb_r|^\g H\left( e, \;\frac{I-e\otimes e}{r \wb_r }+ \left( \frac{ \wb_{rr}  }{ \wb_r^2 } + Z(\wb) \right)e\otimes e\;\right),  \;\;\forall r\ge R. \nonumber
\eea
From parts (ix) of Cases A, B and C of Remark \ref{sec3.15}, we have that
\ben
\frac{r \wb_{rr} }{\wb_r }=\frac{rv^{\prime\prime}(r) }{ v^{\prime}(r) }\le \left\{ \begin{array}{lcr} 1, &\; \;\;\;\s=\g=2,\;\;k=1,\\ \gs-1, & \;\;\;\s=\g>2,\;\;k>1,\\ \gs-1, &\;\;\;\; \s>\g\ge 2,\;\; k\ge 1. \end{array} \right. 
\een
Using the first version in (\ref{sec6.34}) and noting that $\gs\le 2$, $1-rw_{rr}/w_r\ge 0$, one estimates (see (\ref{sec6.180})) 
$$H\left( e, \;\left(\; r|\wb_r|Z(\wb)\;+1   - \frac{r\wb_{rr}}{\wb_r}\right)e\otimes e -I \right)\ge H(e, -I)\ge -|N|,\;\;\;0\le r\le R.$$
Hence, in $0\le r \le R$,
\ben
&&H(D\wb, D^2\wb+Z(\wb) D\wb\otimes D\wb)+\chi(t)|D\wb|^\s-\wb_t \\
&&\qquad\qquad\qquad\qquad \ge - \left[ \frac{ (b(1+T) v^{\prime}(r) )^k|N| }{r}+\al (b(1+T) v^{\prime}(r) )^\s-a- bv(r) \right].
\een
Next, employing the estimate in (\ref{sec5.40}), i.e., $v^{\prime}(r)\le \gs r^{\gs-1}$, and $(\gs-1)k=\gs$, we get, in $0\le r\le R$,
\ben
&&H(D\wb, D^2\wb+Z(\wb) D\wb\otimes D\wb)+\chi(t)|D\wb|^\s-\wb_t \\
&&\qquad\qquad\qquad\qquad \ge - \left[ \frac{ (\gs b(1+T) r^{\gs-1} )^k|N| }{r}+\al (\gs b(1+T)r^{\gs-1} )^\s-a \right]\\
&&\qquad\qquad\qquad\qquad \ge - \left[ (\gs b(1+T) )^k |N|+\al (\gs b(1+T))^\s R^{\s(\gs-1)}-a \right].
\een
As done in (\ref{sec5.100}), we select an appropriate $a$. Thus, $\wb$ is a sub-solution in $B_R(o)\times (0,T)$.

Next, in $r\ge R$, one finds that (see (\ref{sec5.3}))
\bea\label{sec6.37}
&&|\wb_r|^\g H\left( e, \;\frac{I-e\otimes e}{r \wb_r }+ \left( \frac{ \wb_{rr}  }{ \wb_r^2 } + Z(\wb) \right)e\otimes e\;\right)   \nonumber\\
&&\qquad\qquad =( b(1+t)v^{\prime}(r) )^\g H\left( e, \frac{e\otimes e-I}{b(1+t) rv^{\prime}(r) }+\left( Z(\wb)-\frac{ v^{\prime\prime}(r) }{b(1+t) (v^{\prime}(r) )^2 } \right)e\otimes e   \right) \nonumber\\
&&\qquad\qquad \ge b^k( (1+t)v^{\prime}(r) )^\g H\left( e, \left( bZ(\wb)+\frac{1}{rv^{\prime}(r)} -\frac{ v^{\prime\prime}(r) }{(v^{\prime}(r) )^2 } \right)e\otimes e-\frac{I}{rv^{\prime}(r) }   \right),
\eea
where we have factored out $1/b$ and used that $\g=k+1$ and $e\otimes e-I\le 0$. 

We now recall (\ref{sec5.60}) i.e.,
\ben
0<\min\left( \frac{1}{ rv^{\prime}(r) },\;\frac{v^{\prime\prime}(r)}{ (v^{\prime}(r))^2} \right)\le\max\left( \frac{1}{ rv^{\prime}(r) },\;\frac{v^{\prime\prime}(r)}{ (v^{\prime}(r))^2} \right)\le 2,\;\;\mbox{in $r\ge R\ge 1$.}
\een
Employing this estimate in (\ref{sec6.37}) and disregarding the term with $Z$, we get
\ben
|\wb_r|^\g H\left( e, \;\frac{I-e\otimes e}{r \wb_r }+ \left( \frac{ \wb_{rr}  }{ \wb_r^2 } + Z(\wb) \right)e\otimes e\;\right)
\ge b^k( (1+T)v^{\prime}(r) )^\g S,
\een
where 
\ben
S=\min_{|e|=1} H(e, -2(I+e\otimes e) ).
\een
Clearly, by (\ref{sec6.180}), $S\le N<0$ and we get that
\ben
&&H(D\wb, D^2\wb+Z(\wb) D\wb\otimes D\wb)+\chi(t)|D\wb|^\s-\wb_t\\
&&\qquad\qquad\ge - \left[ b^k( (1+T)v^{\prime}(r) )^\g |S|+ \al (b(1+T))^\s (v^{\prime}(r))^\s-a-bv(r) \right],
\een
which is analogous to (\ref{sec5.8}). As done in Section 5, a choice for $b$ (see (\ref{sec5.103})) can now be made. Also,
$\lim_{b\rightarrow 0} a=0$.  $\Box$

\section{Proofs of Theorems \ref{sec2.10}-\ref{sec2.12}}

Let $T>0$ and $(x,t)\in \IR^n\times(0,T),\;n\ge 2$. Set 
$$\mbox{(i)}\;\;\mu=\inf_{\IR^n}h,\quad \nu=\sup_{\IR^n}h,\;\;\mbox{and (ii) assume that $-\infty<\mu\le \nu<\infty$}.$$

Recall that $k=k_1+1,\;\g=k+1$ and $\al=\sup_{[0,T]} |\chi(t)|.$
\vsp
{\bf Proof of Theorem \ref{sec2.10}.} 
Set $r=|x|$ and let $\eta>0$ be small. Choose $\rho>\rho_0$, where $\rho_0$ is large enough so that
\eqRef{sec7.1}
\sup_{ [0,\rho]\times[0,T]}u(x,t)\le \eta \rho^{\bt},\;\;\;\forall\;\rho\ge \rho_0.
\ee
where $\bt$ is as described in the statement of Theorem \ref{sec2.10}.

{\bf Proof of Theorem \ref{sec2.10}(a) $\s=0$:} Recall from (\ref{sec5.0}) that the super-solution $w(x,t)$, with $m=\nu$, is
\eqRef{sec7.2}
w(x,t)=\nu+at+bv(r),\quad \mbox{where $a>0$, $b>0$ and}\quad \lim_{b\rightarrow 0}a=\al.
\ee
For details, see Part I in Section 5, (\ref{sec5.9300}) in Sub-Part (i), (\ref{sec5.970}) in Sub-Part (ii) and (\ref{sec5.98}). Note that
\eqRef{sec7.3}
\mbox{(a)\; If $k=1$ then}\;\; v(r)=\int_0^{r^2} (1+\tau^{\ve/2})^{-1} d\tau,\;\;\mbox{and \;\; (b)\; if $k>1$ then}\;\; v(r)=r^{\gs}.
\ee
See (\ref{sec5.90}) and (\ref{sec5.95}). Also, in (\ref{sec7.1}) and (\ref{sec7.3}), 
\eqRef{sec7.4}
\mbox{(a)\; if $k=1$ then}\;\; \bt=2-\ve,\;\;\mbox{and (b)\; if $k>1$ then}\;\; \bt=\gs=\g/k.
\ee

Recall that $w$, in (\ref{sec7.2}), is a super-solution in $\IR^n_T$ for any $0<b<b_0$, where $b_0$ is small enough, and for an appropriate $a$ that depends on $b$. 

We observe that by part (iv) of Cases A, B and C of Remark \ref{sec3.15},
$v(r)\ge r^\bt /4$, for $r\ge \rho_1$, where $\rho_1$ is large enough. We now choose a fixed
$0<b<b_0$ and take $\eta=b/8$. and let $\rho_0$ stand for the value of $r$ needed for (\ref{sec7.1}) to hold. 

Set $\rho_2=\max(\rho_0, \rho_1)$ and consider a cylinder $B_\rho(o)\times [0,T]$, where $\rho>R_2$. Let $u$ be a sub-solution such that (\ref{sec7.1}) holds. Then $u(x,0)\le h(x)\le \nu,\;\forall x\in \IR^n$. Clearly, $w(x,0)=\nu+b v(r)\ge u(x,0)$, for $|x|\le \rho$. At $|x|=\rho$, we have
$$w(x,t)\ge bv(R)\ge 2\eta  \rho^\bt\ge u(x,t).$$
Thus, $w\ge u$ on the parabolic boundary of $B_{\rho}(o)\times (0,T)$ and Lemma \ref{sec3.18} to conclude that $u(x,t)\le w(x,t)$ in $B_\rho(o)\times(0,T)$ for any large $\rho$, i.e.,
\ben
u(x,t)\le \nu+at +bv(r),\;\;\forall |x|\le \rho.
\een
Letting $\rho\rightarrow \infty$, we see that $u(x,t)\le \nu+at +bv(r)$ in $\IR^n_T$. Since this holds for any small $b$, using (\ref{sec7.2}), we obtain $u(x,t)\le \nu +\al t$. The claim holds.

{\bf Proof of Theorem \ref{sec2.10}(b) $0<\s\le \g$:} The functions $w$, $v$ and $\bt$ are as in (\ref{sec7.2}), (\ref{sec7.3}) and (\ref{sec7.4}). Refer to Part I in Section 5 and see
Sub-Parts (i) and (ii). Arguing as in the proof of Theorem 2.10(a) above, we see that 
$u(x,t)\le \nu+at + b v(r),$ in $\IR^n_T$, for any $b>0$ small enough. Recalling (\ref{sec5.98}), we get that $u(x,t)\le \nu$ and the claim holds.

{\bf Proof of Theorem \ref{sec2.10}(c) $\s>\g$:} Refer to Part II in Section 4. The quantity $\bt=\s/(\s-1)$ in (\ref{sec7.1}). From (\ref{sec5.99})
\ben
w(x,t)=\nu+at+b(1+t) v(r),\;\mbox{where}\; v(r)=\int_0^{ r^{\gs} }  (1+\tau^p)^{-1} d\tau\;\;\mbox{and}\;\;p=\frac{ \s-\g }{  \g(\s-1) },
\een
where $a>0$ and $b>0$. The function $w$ is super-solution in $\IR^n_T$ for any $0<b<b_0$, where $b_0$ is small enough, and an appropriate $a$ that depends on $b$. Moreover, by (\ref{sec5.200}), 
$$\lim_{b\rightarrow 0} a=0.$$
The rest of the proof is similar to the proof of Theorem \ref{sec2.10}(a).   \quad $\Box$
\vsp
\NI{\bf Proof of Theorem \ref{sec2.11}:} 
\vsp
We start with the proofs of parts (a)-(c).

{\bf Proof of Theorem \ref{sec2.11}(a), (b) and (c) $0\le \s\le \g$:} Fix $y\in \IR^n$ and $R>0$ and set $r=|x-y|$.
Recall from (\ref{sec6.23})
$$\bar{w}(x,t)=\mu+v(r)-F t, \;\mbox{where}\;\;v(r)=E\int_0^{\om^2} (\tau^p-1)^{-1} d\tau\;\;\mbox{and}\;\;\om=r/R.$$
See Sub-Cases (a), (b) and (c) of Case I in Section 6. From (\ref{sec6.33}) we see that 
$$\lim_{R\rightarrow \infty}F=\left\{ \begin{array}{ccc}  \al,&\s=0,\\ ( \al^\g / (\ell \cH)^\s )^{1/(\g-\s)},&    0<\s<\g,\\ 0, & \s=\g,\;\al <\ell \cH.  \end{array}\right. $$ 

Let $u$ be as in the theorem. Clearly, $\wb(x,0)\le h(x)\le u(x,0)$ in $\IR^n$. Since 
$\sup|u|<\infty$ in $B_R(y)\times [0,T]$, $\wb(x,t)\le u(x,t)$ on $R^\prime\le |x-y|<R$ for some $R^\prime<R$. By Lemma \ref{sec3.18}, $\wb\le u$ in $B_R(y)\times (0,T).$

Thus, $w(y,t)\le u(y,t)$ and $u(y,t)\ge \mu-Ft$. Letting $R\rightarrow \infty$, we get, 
$$u(y,t)\ge \left\{ \begin{array}{ccc} \mu- \al t, & \s=0,\\ \mu-t(\al^\g/(\ell \cH)^\s)^{1/(\g-\s)},& 0<\s<\g,\\ \mu,& \s=\g,\;\al<\ell \cH. \end{array}\right. $$ 

In order to show the claim for $\chi\ge 0$, take $\al=0$ and refer to Sub-Part (b) in Part I in Section 5. 

{\bf Proof of Theorem \ref{sec2.11}(d) and (e) $\s\ge \g$:} If $\s=\g$ then we take $\al\ge \ell \cH$. 

Assume that 
\eqRef{sec7.5}
\sup_{B_R(o)\times [0,T]} (-u(x,t))\le o(R^\bt),\quad\mbox{as $R\rightarrow \infty.$}
\ee
Recall Sub-Cases (i) and (ii) in Case 2 in Section 6. We take 
$$ \wb(x,t)=\mu-at -b(1+t) v(r),\;\;\mbox{in $\IR^n_T$.}$$

Suppose that $\s=\g$ and $\al\ge \ell \cH$. If (a) $k=1$ and $\g=2$ then $\bt=2-\ve,$ for any small $\ve>0$, in (\ref{sec7.5}), and we take 
$$v(r)=\int_0^{r^2} (1+\tau^{\ve/2})^{-1}\;d\tau,$$
and (b) $k>1$ and $\g>2$ then $\bt=\gs$, in (\ref{sec7.5}), and we take $v(r)=r^{\gs}$. 

If $\s>\g$ and $k\ge 1$ then $\bt=\s/(\s-1)$, in (\ref{sec7.5}), and 
$$v(r)=\int_0^{ r^{\gs}} (1+\tau^{p} )^{-1}\;d \tau \quad \mbox{where}\; \; p=\frac{\s-\g}{\g(\s-1)}.$$
\vsp
It is to be noted that $\lim_{b\rightarrow 0} a=0$ in the situations stated above. The rest of the proof is similar to that of Theorem \ref{sec2.10}.  $\Box$

\vsp
\NI {\bf Proof of Theorem \ref{sec2.12}.} We take $\al=\s=0$ in Theorems \ref{sec2.10}. Let $u>0$ be as in the statement of the theorem. Set
$v=\phi^{-1}(u)$ and $h=\phi^{-1}(g)$, see Lemma \ref{sec3.24}. Let $k\ge 1$. 

{\bf Proof of Theorem \ref{sec2.12}(a):} Since $v$ is a sub-solution we have that $\sup_{B_R(o)\times [0,T]}v(x,t)\le o(R^{\g/k} )$ as $R\rightarrow \infty$. By Lemma \ref{sec3.24},
$v\in usc(\overline{\IR^n_T})$ solves
\ben
&&H(Dv, D^2v+Z(v)Dv\otimes Dv)-v_t\ge 0,\;\mbox{in $\IR^n_T$,}\\
&&\quad\qquad\mbox{and}\;\; v(x,0)\le  \phi^{-1}(h(x)),\;\mbox{for all $x\in \IR^n$,}
\een
By Theorem \ref{sec2.10}(a),
$\sup_{\IR^n_T} v\le \sup_{\IR^n} \phi^{-1}(h),$
and thus, $\sup_{\IR^n_T}u\le \sup_{\IR^n} g.$

{\bf Proof of Theorem \ref{sec2.12}(b):} In this case, $v\in lsc(\overline{\IR^n_T})$ solves
\ben
&&H(Dv, D^2v+Z(v)Dv\otimes Dv)-v_t\ge 0,\;\mbox{in $\IR^n_T$,}\\
&&\quad\qquad\mbox{and}\;\; v(x,0)\le  \phi^{-1}(h(x)),\;\mbox{for all $x\in \IR^n$.}
\een
By Theorem \ref{sec2.11}(a), $\inf_{\IR^n_T} v\ge \inf_{\IR^n} \phi^{-1}(h)$ and hence,
$\inf_{\IR^n_T}u\ge \inf_{\IR^n} h.$ 

The case $k=1$ also follows analogously.
\quad$\Box$

\vsp

\NI Department of Mathematics, Western Kentucky University, Bowling Green, Ky 42101, USA\\
\NI Department of Liberal Arts, Savannah College of Arts and Design, Savannah, GA 31405, USA

\end{document}